\documentclass[MSNbibl,number,citesort,seceqn,dvips]{arxbj}
\usepackage{upgreek,dcolumn}
\usepackage{stfloats}
\usepackage{graphicx}

\aid{0}
\volume{20}
\issue{1}
\pubyear{2014}
\firstpage{109}
\lastpage{140}
\doi{10.3150/12-BEJ478} %kopijuoti is PTS

\makeatletter
\fnbelowfloat
\newcolumntype{d}[1]{D{.}{.}{#1}}

\newcommand{\angler}{\rangle}
\newcommand{\anglel}{\langle}
\renewcommand{\citep}[1]{(\citeauthor{#1} \citeyear{#1})}
\newcommand{\eqref}[1]{(\ref{#1})}
\newcommand{\R}{\mathbb{R}}

\newcommand{\argmin}{\operatorname{argmin}}
\newcommand{\sgn}{\operatorname{sgn}}
\newcommand{\re}{\operatorname{Re}}
\newcommand{\im}{\operatorname{Im}}
\newcommand{\var}{\operatorname{Var}}
\newcommand{\cov}{\operatorname{Cov}}
\newcommand{\supp}{\operatorname{supp}}

\newtheorem{thmm}[lemma]{Theorem}
\newtheorem{lemma}{Lemma}[section]
\newtheorem{prop}[lemma]{Proposition}

\newproclaim{assumption}{Assumption}

\newproclaim{definition}[lemma]{Definition}
\newremark{example}[lemma]{Example}

\newremark{remark}[lemma]{Remark}
\newremark{examples}{Examples}
\makeatother

\begin{document}
\begin{frontmatter}

\title{Calibration of self-decomposable L\'evy models}
\runtitle{Calibration of self-decomposable L\'evy models}

\begin{aug}
%%%% inicialai - be tarpu
\author{\fnms{Mathias} \snm{Trabs}\ead[label=e1]{trabs@math.hu-berlin.de}}
\address{Institut f\"ur Mathematik, Humboldt-Universit\"at zu
Berlin, Unter den Linden 6, 10099 Berlin, Germany. \printead{e1}}
\runauthor{M. Trabs} %% auto
\end{aug}

% HISTORY:
\received{\smonth{11} \syear{2011}}
\revised{\smonth{7} \syear{2012}}

% ABSTRACT
%
\begin{abstract}
We study the nonparametric calibration of exponential L\'evy models
with infinite jump activity.
In particular our analysis applies to self-decomposable processes
whose jump density can
be characterized by the $k$-function, which is typically nonsmooth at
zero. On the one hand the
estimation of the drift, of the activity measure $\alpha:=k(0+)+k(0-)$
and of analogous parameters
for the derivatives of the $k$-function are considered and on the other
hand we estimate nonparametrically
the $k$-function. Minimax convergence rates are derived. Since the rates
depend on $\alpha$, we construct
estimators adapting to this unknown parameter. Our estimation method
is based on spectral representations
of the observed option prices and on a regularization by cutting off
high frequencies. Finally,
the procedure is applied to simulations and real data.
\end{abstract}

% KEYWORDS
% visi is mazosios raides ir pagal abecele
%

\end{frontmatter}

%s1 #&#
\section{Introduction}\label{sec1}
Since Merton \cite{merton1976} introduced his discontinuous asset price
model, stock returns were frequently described by exponentials of L\'
evy processes. A review of recent pricing and hedging results for these
models is given by Tankov \cite{Tankov2011}. The calibration of the
underlying model, that is in the case of L\'evy models the estimation
of the characteristic triplet $(\sigma, \gamma, \nu)$, from historical
asset prices is mostly studied in parametric models only, consider the
survey paper of Eberlein~\cite{eberlein2012} and the references therein.
Remarkable exceptions are the nonparametric penalized least squares
method by Cont and Tankov \cite{contTankov2004} and the spectral calibration
procedure by Belomestny and Rei{\ss}~\cite{reiss12006}. Both articles concentrate on models
of finite jump activity. Our goal is to extend their results to
infinite intensity models. A class which attracted much interest in
financial modeling is given by self-decomposable L\'evy processes,
examples are the hyperbolic model (Eberlein, Keller and
Prause~\cite{eberlein1998}) or the
variance gamma model (Madan and Seneta \cite{madan1990}, Madan, Carr and Chang~\cite{madan1998}). Moreover,
self-decomposable distributions are discussed in the financial
investigation using Sato processes (Carr \textit{et~al.} \cite{carr2007}, Eberlein and Madan \cite{eberlein2009}).
Our results can be applied in this context, too. The nonparametric
calibration of L\'evy models is not only relevant for stock prices, for
instance, it can be used for the Libor market as well (see Belomestny and
Schoenmakers \cite{belomestnySchoenmakers2011}). In the context of Ornstein--Uhlenbeck
processes, the nonparametric inference of self-decomposable L\'evy
processes was considered by Jongbloed, van~der Meulen and van~der
Vaart \cite{Jongbloed2005}.

Owing to the infinite activity, the features of market prices can be
reproduced even without a diffusion part (cf. Carr \textit{et~al.} \cite{carr2002}) and
thus we study pure-jump L\'evy processes. More precisely, we assume
that the jump density satisfies
{\renewcommand{\theequation}{K}
%e1.1 #&#
\begin{equation}
\label{propK} \nu(\mathrm{d}x)=\frac{k(x)}{|x|}\,\mathrm{d}x\qquad
\mbox{where }k\dvtx \mathbb R\to\mathbb R_+\mbox{ has bounded variation.}%\tag{K}
\end{equation}
}
\hspace*{-2pt}When $k$ increases on $(-\infty,0)$ and decreases on $(0,\infty)$, it
is called $k$-function and the processes is self-decomposable. Further
examples which have property \eqref{propK} are compound Poisson
processes and limit distributions of branching processes as considered
by Keller-Ressel and
Mijatovi{\'c} \cite{KellerRessel2012}. Using the bounded variation of $k$, we
show that the estimation problem is only mildly ill-posed. While the
Blumenthal--Getoor index, which was estimated by Belomestny \cite
{belomestny2010}, is zero in our model, the infinite activity can be
described on a finer scale by the parameter
\[
\alpha:=k(0+)+k(0-).
\]
Since $k$ is typically nonsmooth at zero, we face two estimation
problems: First, to give a proper description of $k$ at zero, we
propose estimators for $\alpha$ and its analogs
$k^{(j)}(0+)+k^{(j)}(0-)$, with $j\geq1$, for the derivatives of $k$ as
well as for the drift $\gamma$, which can be estimated similarly. We
prove convergence rates for their mean squared error which turn out to
be optimal in minimax sense up to a logarithmic factor. Second, we
construct a nonparametric estimator of $k$ whose mean integrated
squared error converges with nearly optimal rates.
Owing to bid-ask spreads and other market frictions, we observe only
noisy option prices. The definition of the estimators is based on the
relation between these prices and the characteristic function of the
driving process established by Carr and Madan \cite{carrMadan1999} and on different
spectral representations of the characteristic exponent. Smoothing is
done by cutting off all frequencies higher than a certain value
depending on a maximal permitted parameter $\alpha$. The whole
estimation procedure is computationally efficient and achieves good
results in simulations and in real data examples.
All estimators converge with a polynomial rate, where the maximal
$\alpha$ determines the ill-posedness of the problem. Assuming
sub-Gaussian error distributions, we provide an estimator with $\alpha
$-adaptive rates. The main tool for this result is a concentration
inequality for our estimator $\hat\alpha$ which might be of independent
interest.

This work is organized as follows: In Section \ref{sec2}, we describe the
setting of our estimation procedure and derive the necessary
representations of the characteristic exponent. The estimators are
described in Section \ref{sec3}, where we also determine the convergence rates.
The construction of the $\alpha$-adaptive estimator of $\alpha$ is
contained in Section \ref{secadaptation}. In view of simulations and real data, we
discuss our theoretical results and the implementation of the procedure
in Section \ref{sec5}. All proofs are given in Section \ref{sec6}.

%s2 #&#
\section{The model}\label{sec2}

%s2.1 #&#
\subsection{Self-decomposable L\'evy processes}
A real valued random variable X has a \textit{self-decomposable} law if
for any $b\in(0,1)$ there is an independent random variable $Z_b$ such
that $X\stackrel{d}{=}bX+Z_b$. Since each self-decomposable
distribution is infinitely divisible (see Proposition 15.5 in \cite{sato1999}), we
can define the corresponding \textit{self-decomposable L\'evy process}.
Self-decomposable laws can be understood as the class of limit
distributions of converging scaled sums of independent random variables
(Theorem 15.3 in \cite{sato1999}). This characterization is of economical
interest. If we understand the price of an asset as an aggregate of
small independent influences and release from the $\sqrt{n}$ scaling,
which leads to diffusion models, we automatically end up in a
self-decomposable price process.

Sato \cite{sato1999} shows that the jump
measure of a self-decomposable distribution is always absolutely
continuous with respect to the Lebesgue measure and its density can be
characterized through \eqref{propK} where $k$ needs to be increasing on
$\mathbb R_-$ and decreasing on $\mathbb R_+$. Note that
self-decomposability does not affect the volatility $\sigma$ nor the
drift $\gamma$ of the L\'evy process.

Assuming $\sigma=0$ and property \eqref{propK}, the process $X_t$ has
finite variation and the characteristic function of $X_T$ is given by
the L\'evy--Khintchine representation
%
%
%e2.1 #&#
\begin{equation}
\label{eqphiT} \varphi_T(u):=\mathbb E\bigl[\mathrm{e}^{\mathrm{i}uX_T}
\bigr]=\exp\biggl(T \biggl(\mathrm{i}\gamma u+\int_{-\infty}^\infty
\bigl(\mathrm{e}^{\mathrm{i}ux}-1\bigr)\frac{k(x)}{|x|}\,\mathrm{d}x \biggr)
\biggr).
\end{equation}
Motivated by a martingale argument, we will suppose the exponential
moment condition $\mathbb E[\mathrm{e}^{X_t}]=1$ for all $t\geq0$,
which yields
%
%
%e2.2 #&#
\begin{equation}
\label{eqmartingalCondition} 0=\gamma+\int_{-\infty}^\infty
\bigl(\mathrm{e}^{x}-1\bigr)\frac{k(x)}{|x|}\,\mathrm{d}x.
\end{equation}
In particular, we will impose $\int_{-\infty}^\infty(\mathrm
{e}^x-1)\frac
{k(x)}{|x|}\,\mathrm{d}x<\infty$. In this case, $\varphi_T$ is
defined on the
strip $\{z\in\mathbb C|\im z\in[-1,0]\}$.

Besides L\'evy processes there is another class that is closely related
to self-decompos\-ability. Assuming self-similarity, that means
$(Y_{at})\stackrel{d}{=}(a^HY_{t})$, for all $a>0$ and some exponent
$H>0$, instead of stationary increments, $Y_t$ is a \textit{Sato
processes}. Sato \cite{sato1991} showed that self-decomposable
distributions can be characterized as the laws at unit time of these
processes. From the self-similarity and self-decomposability follows
for $T>0$
\[
\varphi_{Y_T}(u)=\mathbb E\bigl[\mathrm{e}^{\mathrm{i}uY_T}\bigr]=
\mathbb E\bigl[\mathrm{e}^{\mathrm{i}T^HuY_1}\bigr]=\exp\biggl(\mathrm{i}T^H\gamma
u+\int_{-\infty}^\infty\bigl(\mathrm{e}^{\mathrm{i}ux}-1
\bigr)\frac
{k(T^{-H}x)}{|x|}\,\mathrm{d} x \biggr).
\]
Since our estimation procedure only depends through equation \eqref
{eqphiT} on the distributional structure of the underlying process, we
can apply the estimators directly to Sato processes using
$T_{s}=1,\gamma_{s}=T^H\gamma$ and $k_{s}(\cdot)=k(T^{-H}\bullet)$
instead of $T$, $\gamma$ and $k$. However, we concentrate on L\'evy
processes in the sequel.

For self-decomposable distributions the parameter $\alpha$ captures
many of its properties such as the smoothness of the densities of the
marginal distributions (Theorem 28.4 in~\cite{sato1999}) and the tail
behavior of the characteristic function. This holds even for the more
general class of L\'evy processes that satisfy property \eqref{propK}.
Recall that $k$ has bounded variation if and only if
\[
\|k\|_{\mathrm{TV}}:=\sup\Biggl\{\sum_{i=1}^n\bigl|k(x_i)-k(x_{i-1})\bigr|\dvt n
\in\mathbb N,-\infty<x_0<\cdots<x_n<\infty\Biggr\}<
\infty.
\]
In particular, $\|k\|_{\mathrm{TV}}<\infty$ implies $\alpha<\infty$. Similarly
to deconvolution problems, the stochastic error in our model is driven
by $|\varphi_T(u-\mathrm{i})|^{-1}$ and thus we prove the following lemma in the
\hyperref[app]{Appendix}.
%
%
%le2.1 #&#
\begin{lemma}\label{lemboundPhi}
Let $X_t$ have property \textup{\eqref{propK}} and $\sigma=0$ and let the
martingale condition \eqref{eqmartingalCondition} hold.
\begin{longlist}[(ii)]
\item[(i)] If $\|\mathrm{e}^xk(x)\|_{L^1}<\infty$ and $q_k:=\sup_{x\in(0,1]}\frac
{k(x)+k(-x)-\alpha}{x}<\infty$ then there exists a constant
$C_\varphi
=C_\varphi(T,\max\{q_k,\|\mathrm{e}^xk(x)\|_{L^1},\|k\|_{\mathrm{TV}}\})>0$ such that for
all $u\in\mathbb R$ with $|u|\geq1$ we obtain the bound
\[
\bigl|\varphi_T(u-\mathrm{i})\bigr|\geq C_\varphi|u|^{-T\alpha}.
\]
\item[(ii)] Let $\bar\alpha,R>0$ then $|\varphi_T(u-\mathrm{i})|\geq C_\varphi
(T,R)|u|^{-T\bar\alpha}$ holds uniformly over all $|u|\ge1$ and all
$X_T$ with $\alpha\leq\bar\alpha$ and $\max\{q_k,\|\mathrm{e}^xk(x)\|_{L^1},\|k\|
_{\mathrm{TV}}\}\leq R$.
\end{longlist}
\end{lemma}
The value $q_k$ as defined in the lemma can be understood as the
largest slop of $k$ near zero. If the process is self-decomposable than
$q_k\le0$ holds and the bounded variation norm equals~$\alpha$.
Otherwise, we can use $q_k\le\sup_{|x|\le1}|k'(x)|$ and $\|k\|_{\mathrm{TV}}\le\|
k'\|_{L^1}$, assuming the derivative $k'$ exists, is bounded on
$[-1,1]$ and integrable on $\mathbb R$. If either $\sigma>0$ or
property \eqref{propK} is violated, $\varphi_T$ can decay faster than
any polynomial order, for example, consider self-decomposable processes
with $\alpha=\infty$ (see \cite{sato1999}, Lemma 28.5). Hence, the conditions
of Lemma \ref{lemboundPhi} are sharp.

%s2.2 #&#
\subsection{Asset prices and Vanilla options}
Let $r\geq0$ be the risk-less interest rate in the market and $S_0>0$
denote the initial value of the asset. In an exponential L\'evy model
the price process is given by
\[
S_t=S_0\mathrm{e}^{rt+X_t},
\]
where $X_t$ is a L\'evy process described by the characteristic triplet
$(\sigma,\gamma,\nu)$. Throughout these notes, we assume $X_t$ has
property \eqref{propK} and $\sigma=0$. On the probability space
$(\Omega
,\mathcal F,\mathbb P)$ with pricing (or martingale) measure $\mathbb
P$ the discounted process $(\mathrm{e}^{-rt}S_t)$ is a martingale with respect
to its natural filtration $(\mathcal F_t)$. This is equivalent to
$\mathbb E[\mathrm{e}^{X_t}]=1$ for all $t\geq0$ and thus, the martingale
condition \eqref{eqmartingalCondition} holds.

At time $t=0$ the risk neutral price of an European call option with
underlying $S$, time to maturity $T$ and strike price $K$ is given by
$\mathrm{e}^{-rT}\mathbb E[(S_T-K)_+],$ where $A_+:=\max\{0,A\}$, and similarly
$\mathrm{e}^{-rT}\mathbb E[(K-S_T)_+]$ is the price of European put. In terms of
the negative log-forward moneyness $x:=\log(K/S_0)-rT$ the prices can
be expressed as
\[
\mathcal C(x,T)=S_0\mathbb E\bigl[\bigl(\mathrm{e}^{X_T}-\mathrm{e}^x
\bigr)_+\bigr] \quad\mbox{and}\quad \mathcal P(x,T)=S_0\mathbb E\bigl[
\bigl(\mathrm{e}^x-\mathrm{e}^{X_T}\bigr)_+\bigr].
\]
Carr and Madan \cite{carrMadan1999} introduced the option function
\[
\mathcal{O}(x):=\cases{ %
 S_0^{-1}
\mathcal{C}(x,T), &\quad  $x\ge0,$
\vspace*{2pt}\cr
S_0^{-1}\mathcal{P}(x,T), &\quad  $x<0$,}
\]
and set the Fourier transform $\mathcal{FO}(u):=\int_{-\infty
}^\infty
\mathrm{e}^{\mathrm{i}ux}\mathcal O(x)\,\mathrm{d}x$ in relation to the characteristic function
$\varphi_T$ through the pricing formula
%
%
%e2.3 #&#
\begin{equation}
\label{eqpricingFormula} \mathcal{FO}(u)=\frac{1-\varphi
_T(u-\mathrm{i})}{u(u-\mathrm{i})}, \qquad u\in\mathbb R
\setminus\{0\}.
\end{equation}
The properties of $\mathcal O$ were studied further by Belomestny and
Rei{\ss} \cite{reiss12006}.
%: At any $x\in\mathbb R\setminus\{0\}$ the function $\mathcal O$ is
%twice differentiable with $\int_{\mathbb R}|\mathcal O''(x)|\di x
%at zero. Additionally, they showed that Assumption \ref{assmoment}
%ensures an exponential decay of the option function, i.e. $|\mathcal
%O(x)|\lesssim e^{-|x|}$ holds for $x\in\mathbb R$.
In particular, they showed that the option function is contained in
$C^1(\mathbb R\setminus\{0\})$ and decays exponentially under
the following assumption.
%
%
%as1 #&#
\begin{assumption}\label{assmoment}
We assume that $C_2:=\mathbb E[\mathrm{e}^{2X_T}]$ is finite, which is
equivalent to the moment condition $\mathbb E[S_t^2]<\infty$.
\end{assumption}
Our observations are given by
%
%
%e2.4 #&#
\begin{equation}
\label{eqregressionModel} O_j=\mathcal{O}(x_j)+
\delta_j\varepsilon_j,\qquad  j=1,\ldots,N,
\end{equation}
where the noise $(\varepsilon_j)$ consists of independent, centered
random variables with $\mathbb E[\varepsilon_j^2]=1$ and $\sup_j\mathbb
E[\varepsilon_j^4]<\infty$. The noise levels $\delta_j$ are assumed to
be positive and known. In practice, the uncertainty is due to market
frictions such as bid-ask spreads.

%s2.3 #&#
\subsection{Representation of the characteristic exponent}
Using \eqref{eqphiT} and \eqref{eqpricingFormula}, the shifted
characteristic exponent is given by
%
%
%e2.5 #&#
%e2.6 #&#
\begin{eqnarray}
\psi(u)&:=&\frac{1}{T}\log\bigl(1+\mathrm{i}u(1+\mathrm{i}u)\mathcal{FO}(u)\bigr)=
\frac
{1}{T}\log\bigl(\varphi_T(u-\mathrm{i})\bigr)\label{eqdefPsi}
\\
& =& \mathrm{i}\gamma u + \gamma+\int_{-\infty}^\infty
\bigl(\mathrm{e}^{\mathrm{i}(u-\mathrm{i})x}-1\bigr)\frac
{k(x)}{|x|}\,\mathrm{d}x\label{eqrepPsi0}
\end{eqnarray}
for $u\in\mathbb R$. Note that the last line equals zero for $u=0$
because of the martingale condition \eqref{eqmartingalCondition}.
Throughout, we choose a distinguished logarithm, that is a version of
the complex logarithm such that $\psi$ is continuous with $\psi(0)=0$.
Under the assumption that $\int_{-\infty}^{\infty}(1\vee \mathrm{e}^x)k(x)\,
\mathrm{d}
x$\footnote{We denote $A\wedge B:=\min\{A,B\}$ and $A\vee B:=\max\{
A,B\}
$ for $A,B\in\mathbb R$.} is finite, we can apply Fubini's theorem to obtain
%
%
%e2.7 #&#
\begin{equation}
\psi(u)=\mathrm{i}\gamma u + \gamma+\int_0^1 \mathrm{i}(u-\mathrm{i})
\mathcal{F}\bigl(\sgn(x)k(x)\bigr) \bigl((u-\mathrm{i})t\bigr)\,\mathrm{d}t,
\label{eqrepPsi1}
\end{equation}
where the Fourier transform $\mathcal F(\sgn\cdot k)$ is well defined
on $\{z\in\mathbb C|\im z\in[-1,0]\}$. Typically, the~$k$ and its
derivatives are not continuous at zero. Moreover, if $\alpha\neq0$ the
function $x\mapsto\sgn(x)k(x)$ has a jump at zero in every case.
Therefore, the Fourier transform decreases very slowly.
Let $k$ be smooth on $\mathbb R\setminus\{0\}$ and fulfill an
integrability condition which will be important later:
%
%
%as2 #&#
\begin{assumption}\label{assk}
Assume $k\in C^s(\mathbb{R}\setminus\{0\})$ with all derivatives having
a finite right- and left-hand limit at zero and $(1\vee \mathrm{e}^x)k(x), \ldots
, (1\vee \mathrm{e}^x)k^{(s)}(x)\in L^1(\mathbb R)$.
\end{assumption}
To compensate those discontinuities, we add a linear combination of the
functions $h_j(x):=x^je^{-x}\mathbf1_{[0,\infty)}(x), x\in\mathbb R$,
for $j=\mathbb N\cup\{0\}$. Since $h_j\in C^{j-1}(\mathbb{R})$ for
$j\ge
1$ and all $h_j$ are smooth on $\mathbb{R}\setminus\{0\}$, we can find
$\alpha_j, j=0,\ldots,s-2$, such that $\sgn(x)k(x)-\sum_{j=0}^{s-2}\alpha
_j h_j(x)$ is contained in $C^{s-2}(\mathbb{R})\cap C^s(\mathbb
{R}\setminus\{0\})$. This approach yields the following representation.
The proof is given in the supplementary article \cite{trabs2011Supplement}.
%
%
%pr2.2 #&#
\begin{prop}\label{proprepPsi}
Let $s\geq2$. On Assumption \ref{assk}, there exist functions $D\dvtx \{
-1,1\}\to\mathbb C$ and $\rho\dvtx \mathbb R\setminus\{0\}\to\mathbb C$ such
that $|u^{s-1}\rho(u)|$ is bounded in $u$ and it holds
%
%
%e2.8 #&#
\begin{equation}
\psi(u)=D\bigl(\sgn(u)\bigr)+\mathrm{i}\gamma u-\alpha_0\log\bigl(|u|\bigr)+\sum
_{j=1}^{s-2}\frac{\mathrm{i}^j(j-1)!\alpha_j}{u^j}+\rho(u),\qquad u
\neq0.\label{eqrepresentationPsi}
\end{equation}
The coefficients are given by $\alpha_j=\frac{1}{j!}
(k^{(j)}(0+)+k^{(j)}(0-) )-\sum_{m=1}^j\frac{(-1)^{m}}{m!}\alpha_{j-m},$
especially $\alpha_0=\alpha$ holds.
\end{prop}
Representation \eqref{eqrepresentationPsi} allows us to estimate
$\gamma$ and $\alpha_0,\ldots,\alpha_{s-2}$. A plug-in approach yields
estimators for $k^{(j)}(0+)+k^{(j)}(0-) , j=0,\ldots,s-2$. Since we only
apply this representation when $\psi$ is multiplied with weight
functions having roots of degree $s-1$ at zero, the poles that appear
in \eqref{eqrepresentationPsi} do no harm.

Proposition \ref{proprepPsi} covers the case $s\geq2$. For $s=1,$ we
conclude from \eqref{eqrepPsi0}, the martingale condition \eqref
{eqmartingalCondition} and Assumption \ref{assk}
%
%
%e2.9 #&#
\begin{eqnarray}\label{eqrepPsi2}
\psi(u)=\mathrm{i}\gamma u+\int_{-\infty}^\infty
\bigl(\mathrm{e}^{\mathrm{i}ux}-1\bigr)\mathrm{e}^x\frac
{k(x)}{|x|}\,\mathrm{d}x =\mathrm{i}
\gamma u+\mathrm{i}\int_0^u\mathcal F \bigl(
\sgn(x)\mathrm{e}^xk(x) \bigr) (v)\,\mathrm{d}v.
\end{eqnarray}
Hence, $\psi$ is a sum of a constant from the integration, the linear
drift $\mathrm{i}\gamma u$ and a remainder of order $\log|u|$, which follows
from the decay of the Fourier transform as $|u|^{-1}$.
Corollary 8
in \cite{trabs2011Supplement} even shows, that there exists no
$L^2$-consistent estimator of $\alpha$ for $s=1$. Therefore, we
concentrate on the case $s\geq2$ in the sequel.

Equation \eqref{eqrepPsi2} allows another useful observation. Defining
the exponentially scaled $k$-function
\[
k_e(x):=\sgn(x)\mathrm{e}^xk(x), \qquad x\in\mathbb R,
\]
we obtain by differentiation
%
%
%e2.10 #&#
\begin{equation}
\label{eqrepPsiPrim} \psi'(u) =\frac{1}{T}\frac{(\mathrm{i}-2u)\mathcal
{FO}(u)-(u+\mathrm{i}u^2)\mathcal F
(x\mathcal O(x) )(u)}{1+(\mathrm{i}u-u^2)\mathcal{FO}(u)}
=\mathrm{i}\gamma+\mathrm{i}\mathcal Fk_e(u).
\end{equation}
Using this relation, we can define an estimator of $k_e$.

%s3 #&#
\section{Estimation procedure}\label{sec3}
%s3.1 #&#
\subsection{Definition of the estimators and weight functions}\label
{subsecweightFunctions}
Given the observations $\{(x_1,O_1),\ldots,(x_N,O_N)\}$, we fit a
function $\tilde{\mathcal{O}}$ to these data using linear $B$-splines
\[
b_j(x):=\frac{x-x_{j-1}}{x_j-x_{j-1}}\mathbf1_{[x_{j-1},x_j)}+
\frac
{x_{j+1}-x}{x_{j+1}-x_j}\mathbf1_{[x_j,x_{j+1}]}, \qquad j=1,\ldots,N,
\]
and a function $\beta_0$ with $\beta_0'(0+)-\beta_0'(0-)=-1$
to take care of the jump of $\mathcal O'$:
\[
\tilde\mathcal{O}(x)=\beta_0(x)+\sum
_{j=1}^NO_jb_j(x),\qquad  x\in
\mathbb R.
\]
We choose $\beta_0$ with support $[x_{j_0-1},x_{j_0}]$ where $j_0$
satisfies $x_{j_0-1}<0\leq x_{j_0}$. Replacing $\mathcal O$ with~$\tilde\mathcal{O}$ in the representations \eqref{eqdefPsi} and
\eqref
{eqrepPsiPrim} of $\psi$ and $\psi'$, respectively, allows us to
define their empirical versions through
\begin{eqnarray*}
\tilde\psi(u)&:=&\frac{1}{T}\log\bigl(v_{\kappa
(u)}\bigl(1+\mathrm{i}u(1+\mathrm{i}u)
\mathcal{F}\tilde\mathcal{O}(u)\bigr) \bigr),
\\
\tilde{\psi'}(u)&:=&\frac{1}{T}\frac{(\mathrm{i}-2u)\mathcal{F}\tilde\mathcal{O}(u)-
(u+\mathrm{i}u^2)\mathcal F (x\tilde\mathcal{O}(x)
)(u)}{v_{\kappa
(u)}(1+\mathrm{i}u(1+\mathrm{i}u)\mathcal{F}\tilde\mathcal{O}(u))},\qquad  u\in
\mathbb{R},
\end{eqnarray*}
where $\kappa$ is a positive function and we apply a trimming function
given by
\[
v_{\kappa}(z)\dvtx \mathbb C\setminus\{0\}\to\mathbb C,\qquad  z\mapsto\cases{ z,&\quad
$|z|\geq\kappa,$\vspace*{2pt}
\cr
\kappa z/|z|,&\quad  $|z|<\kappa$,}
\]
to stabilize for large stochastic errors. A reasonable choice of
$\kappa
$ will be derived below. The function $\tilde\psi$ is well defined on
the interval $[-U,U]$ on the event
\[
A:=\bigl\{\omega\in\Omega\dvt  1+\mathrm{i}u(1+\mathrm{i}u)\mathcal{F}\bigl(\tilde\mathcal{O}(\omega
,\bullet)
\bigr) (u)\neq0\ \forall u\in[-U,U]\bigr\}\in\mathcal F.
\]
For $\omega\in\Omega\setminus A$, we set $\tilde\psi$ arbitrarily, for
instance equal to zero. The more $\tilde\mathcal{O}$ concentrates
around the true function $\mathcal O$ the greater is the probability of
$A$. S{\"o}hl \cite{soehl2010} shows even that in the continuous-time L\'evy
model with finite jump activity the identity $\mathbb P(A)=1$
holds.

In the spirit of Belomestny and Rei{\ss} \cite{reiss12006}, we estimate the parameters
$\gamma
$ and $\alpha_j, j=0,\ldots,  s-2$, as coefficients of the different
powers of $u$ in equation \eqref{eqrepresentationPsi}. Using a
spectral cut-off value $U>0$, we define
\[
\hat\gamma:=\int_{-U}^U\im\bigl(\tilde\psi(u)
\bigr)w_\gamma^U(u)\,\mathrm{d}u
\]
and for $0
\leq j\leq s-2$
\[
\hat\alpha_j:=\cases{\displaystyle \int_{-U}^U
\re\bigl(\tilde\psi(u)\bigr)w_{\alpha_j}^U(u)\,\mathrm{d}u,&\quad
$\mbox{if }j\mbox{ is even},$\vspace*{2pt}
\cr
\displaystyle\int
_{-U}^U\im\bigl(\tilde\psi(u)
\bigr)w_{\alpha_j}^U(u)\,\mathrm{d}u,&\quad $\mbox{otherwise}.$ }
\]

The weight functions $w_\gamma^U$ and
$w_{\alpha_j}^U$ are chosen such that they filter the
coefficients of interest. Owing to \eqref{eqrepPsiPrim}, the nonparametric
object $k_e$ can be estimated by
%e3.1 #&#
\begin{equation}
\hat k_e(x):=\cases{ \mathcal F^{-1} \bigl[\bigl(-\hat
\gamma-\mathrm{i}\tilde\psi'(u)\bigr)\mathcal FW_k(u/U)
\bigr](x), &\quad $x>0$,\vspace*{2pt}
\cr
\mathcal F^{-1} \bigl[\bigl(-\hat
\gamma-\mathrm{i}\tilde\psi'(u)\bigr)\mathcal FW_k(-u/U)
\bigr](x), &\quad $x<0,$ }\label{eqDefKe}
\end{equation}
applying a one-sided kernel function $W_k$ with bandwidth $U^{-1}$
since we know that $k_e$ jumps only at zero. The condition on the
weights are summarized in the following:
%
%
%as3 #&#
\begin{assumption}\label{assweightFunctions}
We assume:
\begin{itemize}
\item$w_\gamma^U$ fulfills for all odd $j\in\{1,\ldots,s-2\}$
\begin{eqnarray*}
\int_{-U}^Uuw_\gamma^U(u)
\,\mathrm{d}u=1, \qquad\int_{-U}^Uu^{-j}w_\gamma^U(u)
\,\mathrm{d}u=0 \quad\mbox{and} \quad\int_0^Uw_\gamma^U(
\pm u)\,\mathrm{d}u=0.
\end{eqnarray*}
\item$w_{\alpha_0}^U$ satisfies for all even $j\in\{1,\ldots,s-2\}$
\[
\int_{-U}^U\log\bigl(|u|\bigr)w_{\alpha_0}^U(u)
\,\mathrm{d}u=-1,\qquad \int_{-U}^Uu^{-j}w_{\alpha_0}^U(u)
\,\mathrm{d}u=0 \quad\mbox{and}\quad \int_0^Uw_{\alpha_0}^U(
\pm u)\,\mathrm{d}u=0.
\]
\item For $j=1,\ldots,s-2$ the weight functions $w_{\alpha_j}^U$
fulfill\footnote{For $a\in\mathbb R$ let $\lfloor a\rfloor$ denote the
largest integer which is smaller than $a$.}
\begin{eqnarray*}
\int_{-U}^Uu^{-j}w_{\alpha_j}^U(u)
\,\mathrm{d}u&=&\frac{(-1)^{\lfloor
j/2\rfloor
}}{(j-1)!},\qquad \int_{-U}^Uu^{-l}w_{\alpha_j}^U(u)
\,\mathrm{d}u=0 \quad\mbox{and}\\
 \int_0^Uw_{\alpha_j}^U(
\pm u)\,\mathrm{d}u&=&0,
\end{eqnarray*}
where $1\leq l\leq s-2$ and $l$ is even for even $j$ and odd otherwise.
For even $j$ we impose additionally
\[
\int_{-U}^U\log\bigl(|u|\bigr)w_{\alpha_j}^U(u)
\,\mathrm{d}u=0.
\]
\item$W_k$ is of Sobolev smoothness $T\bar\alpha+2$, that is, $\int
(1+|u|^2)^{T\bar\alpha+2}|\mathcal FW_k(u/U)|^2\,\mathrm{d}u<\infty
$, has
support $\supp W_k\subseteq(-\infty,0]$ and fulfills for $l=1,\ldots,s-1$
\[
\int_{\mathbb R} W_k(x)\,\mathrm{d}x=1, \qquad\int
_{\mathbb R} x^lW_k(x)\, \mathrm{d} x=0\quad
\mbox{and}\quad x^{2s-1}W_k(x)\in L^1(\mathbb R).
\]
\end{itemize}
Furthermore, we assume continuity and boundedness of the functions
$u\mapsto u^{-s+1}w_q^1(u)$ for $q\in\{\gamma,\alpha_0,\ldots,\alpha
_{s-2}\}$.
\end{assumption}
The integral conditions can be provided by rescaling: Let $w_q^1$
satisfy Assumption \ref{assweightFunctions} for $q\in\{\gamma
,\alpha_0,\ldots,\alpha_{s-2}\}$ and $U=1$. Since $1=\int_{-1}^1uw_\gamma
^1(u)\,\mathrm{d}u=\int_{-U}^UuU^{-2}w_\gamma^1(u/U)\,\mathrm{d}u$,
we can choose
$w_\gamma^U(u):=U^{-2}w_\gamma^1(\frac{u}{U})$. Similarly, a rescaling
is possible for $w_{\alpha_0}^U$:
\begin{eqnarray*}
-1&=&\int_{-1}^1\log\bigl(|u|\bigr)w_{\alpha_0}^1(u)
\,\mathrm{d}u=\int_{-U}^U\log
\bigl(|u|\bigr)U^{-1}w_{\alpha_0}^1\biggl(\frac{u}{U}
\biggr)\,\mathrm{d}u-\frac{\log
(U)}{U}\int_{-U}^Uw_{\alpha_0}^1
\biggl(\frac{u}{U}\biggr)\,\mathrm{d}u
\\
&=&\int_{-U}^U\log\bigl(|u|\bigr)U^{-1}w_{\alpha_0}^1
\biggl(\frac{u}{U}\biggr)\,\mathrm{d}u.
\end{eqnarray*}
Therefore, we define $w_{\alpha_0}^U(u):=U^{-1}w_{\alpha_0}^1(\frac
{u}{U})$ and analogously $w_{\alpha_j}^U(u):=U^{j-1}w_{\alpha
_j}^1(\frac
{u}{U})$.
The continuity condition on $w^1_q$ in Assumption \ref
{assweightFunctions} is set to take advantage of the decay of the
remainder $\rho$. In combination with the rescaling it implies
%
%
%e3.2 #&#
\begin{equation}
\label{eqweightfunctionProperty} \bigl|w_\gamma^U(u)\bigr|\lesssim
U^{-s-1}|u|^{s-1} \quad\mbox{and}\quad \bigl|w_{\alpha_j}^U(u)\bigr|
\lesssim U^{-s+j}|u|^{s-1},\qquad j=0,\ldots,s-2.
\end{equation}
Throughout, we write $A\lesssim B$ if there is a constant $C>0$
independent of all parameters involved such that $A\le CB$. In the
sequel we assume that the weight functions satisfy Assumption \ref
{assweightFunctions} and the property \eqref
{eqweightfunctionProperty}.

We reduce the loss of $\hat k_e$ by truncating positive values on
$\mathbb R_-$ and negative ones on $\mathbb R_+$. In the
self-decomposable framework there are additional shape restrictions of
the $k$-function which the proposed\vadjust{\goodbreak} estimator $\hat k_e$ does not take
into account. The monotonicity can be generated by a rearrangement of
the function. To this end let $\hat k(x):=(\sgn(x)\mathrm{e}^{-x}\hat
k_e(x)\vee
0)\mathbf1_{[-C,C]}(x), x\in\mathbb R$, where we bounded the support
with an arbitrary large constant $C>0$. The rearranged estimator which
is increasing on $\mathbb R_-$ and decreasing on $\mathbb R_+$ is then
given by
%
%
%e3.3 #&#
\begin{equation}
\label{eqRearranged} \hat k^*(x):=\cases{\displaystyle \inf\biggl\{y\in\mathbb R_+
\Big|\int
_0^C\mathbf1_{\{\hat
k(z)\ge y\}}\,\mathrm{d}z\le
x \biggr\}, &\quad $x>0,$\vspace*{2pt}
\cr
\displaystyle\inf\biggl\{y\in\mathbb R_+ \Big|\int
_0^C\mathbf1_{\{\hat
k(-z)\ge y\}}\,\mathrm{d}z
\le|x| \biggr\}, &\quad $x<0.$ }
\end{equation}
Chernozhukov, Fern{\'a}ndez-Val and
Galichon \cite{Chernozhukov2009} show that the rearrangement reduces weakly
the error for increasing target functions on compact subsets. This
result carries over to our estimation problem.

%s3.2 #&#
\subsection{Convergence rates}
To ensure a well-defined procedure, an exponential decay of $\mathcal
O$, the identity \eqref{eqrepPsi2} and to obtain a lower bound of
$|\varphi_T(u-\mathrm{i})|$, we consider the class $\mathcal G_0(R,\bar\alpha)$.
Uniform convergence results for the parameters will be derived in the
smoothness class $\mathcal G_s(R,\bar\alpha)$.
%
%
%de3.1 #&#
\begin{definition} Let $s\in\mathbb N$ and $R,\bar\alpha>0$. We define
\begin{longlist}[(ii)]
\item[(i)]$\mathcal G_0(R,\bar\alpha)$ as the set of all pairs $\mathcal
P=(\gamma,k)$ where $k$ is of bounded variation and the corresponding
L\'evy process $X$ given by the triplet $(0, \gamma, k(x)/|x|)$
satisfies Assumption~\ref{assmoment} with $C_2\leq R$, martingale
condition \eqref{eqmartingalCondition} as well as
\[
\alpha\in[0,\bar\alpha] \quad\mbox{and}\quad \max\biggl\{\sup_{x\in(0,1]}\biggl\{
\frac{k(x)+k(-x)-\alpha}{x}\biggr\},\bigl\|k_e(x)\bigr\|_{L^1},\|k
\|_{\mathrm{TV}} \biggr\}\leq R,
\]
\item[(ii)]$\mathcal G_s(R,\bar\alpha)$ as the set of all pairs $\mathcal
P=(\gamma,k)\in\mathcal G_0(R,\bar\alpha)$ satisfying additionally
Assumption~\ref{assk} with
\begin{eqnarray*}
\bigl|k^{(l)}(0+)+k^{(l)}(0-)\bigr|&\leq& R,\qquad \mbox{for }l=1,\ldots,s-1,
\\
\bigl\|\bigl(1\vee \mathrm{e}^x\bigr)k^{(l)}(x)\bigr\|_{L^1}&\leq& R,\qquad
\mbox{for }l=0,\ldots,s.
\end{eqnarray*}
\end{longlist}
\end{definition}
In the class $\mathcal G_0(R,\bar\alpha)$ Lemma \ref{lemboundPhi}(ii)
provides a common\vspace*{1pt} lower bound of $|\varphi_T(u-\mathrm{i})|$ for $|u|\geq1$.
Using $\max_{x\in\mathbb R}\frac{1-\cos(x)}{x}\in(0,1]$, we estimate
roughly for $u\in(-1,1)\setminus\{0\}$:
\begin{eqnarray*}
\bigl|\varphi_T(u-\mathrm{i})\bigr|=\exp\biggl(T\int_{-\infty}^\infty
\frac{\cos
(x)-1}{x}\mathrm{e}^{x/|u|}k\bigl(x/|u|\bigr)\,\mathrm{d}x \biggr)\geq
\mathrm{e}^{-TR}.
\end{eqnarray*}
Hence, the choice
\[
\kappa(u):=\kappa_{\bar\alpha}(u):=\cases{ \displaystyle\tfrac{1}{3}\mathrm{e}^{-TR},&\quad
$|u|<1,$\vspace*{2pt}
\cr
\displaystyle\tfrac{1}{3}C_\varphi(T,R)|u|^{-T\bar\alpha},&\quad
$|u|\geq1,$ }
\]
satisfies
%
%
%e3.4 #&#
\begin{equation}
\label{eqconditionKappa} \tfrac{1}{3}\bigl|\varphi_T(u-\mathrm{i})\bigr|\geq
\kappa(u),\qquad  u\in\mathbb R,
\end{equation}
where the factor $1/3$ is used for technical reasons. As discussed
above, we can restrict our investigation to the case $s\geq2$. Since
the L\'evy process is only identifiable if $\mathcal O$ is known on the
whole real line, we consider asymptotics of a growing number of
observations with
\[
\Delta:=\max_{j=2,\ldots,N}(x_{j}-x_{j-1})\to0 \quad\mbox{and}\quad
A:=\min(x_N,-x_1)\to\infty.
\]
Taking into account the numerical interpolation error and the
stochastic error, we analyze the risk of the estimators in terms of the
abstract noise level
\[
\varepsilon:=\Delta^{3/2}+\Delta^{1/2}\|\delta
\|_{l^\infty}.
\]

%th3.2 #&#
\begin{thmm}\label{thmupperBoundsParameters}
Let $s\geq2, R,\bar\alpha>0$ and assume $\mathrm{e}^{-A}\lesssim\Delta^2$ and
$\Delta\|\delta\|_{l^2}^2\lesssim\|\delta\|_{l^\infty}^2$. We choose
the cut-off value $U_{\bar\alpha}:=\varepsilon^{-2/(2s+2T\bar\alpha
+1)}$ to obtain the uniform convergence rates
\begin{eqnarray*}
\sup_{\mathcal P=(\gamma,k)\in\mathcal G_s(R,\bar\alpha
)}\mathbb E_{\mathcal P}\bigl[|\hat\gamma-
\gamma|^2\bigr]^{1/2} &\lesssim&\varepsilon^{2s/(2s+2T\bar\alpha+1)}
\quad\mbox{and}
\\
\sup_{\mathcal P=(\gamma,k)\in\mathcal G_s(R,\bar\alpha
)}\mathbb E_{\mathcal P}\bigl[|\hat\alpha_j-
\alpha_j|^2\bigr]^{1/2} &\lesssim&
\varepsilon^{2(s-1-j)/(2s+2T\bar\alpha+1)},\qquad j=0,\ldots,s-2.
\end{eqnarray*}
\end{thmm}
As one may expect the rates for $\alpha_j,j=0,\ldots,s-2,$ become slower
as $j$ gets closer to its maximal value because the profit from the
smoothness of $k$ decreases. Note that the cut-off for all estimators
is the same. In contrast to $\mathcal G_s(R,\bar\alpha)$ we assume
Sobolev conditions on $k_e$ in the class $\mathcal H_s(R,\bar\alpha)$
in order to apply $L^2$-Fourier analysis.
%
%
%de3.3 #&#
\begin{definition}
Let $s\in\mathbb N$ and $R,\bar\alpha>0$. We define $\mathcal
H_s(R,\bar\alpha)$ as the set of all pairs $\mathcal P=(\gamma,k)\in
\mathcal G_0(R,\bar\alpha)$ satisfying additionally $k\in C^s(\mathbb
R\setminus\{0\})$, $\mathbb E_{\mathcal P}[|X_Te^{X_T}|]\leq R$ for
corresponding L\'evy process $X$ as well as
\begin{eqnarray*}
|\gamma|\leq R\quad\mbox{and} \quad\bigl\|k_e^{(l)}\bigr\|_{L^2}
\leq R\qquad \mbox{for }l=0,\ldots,s.
\end{eqnarray*}
\end{definition}
In the next theorem the conditions on $A$ and $\delta$ are stronger
than for the upper bounds of the parameters which is due to the
necessity to estimate also the derivative of $\psi$. However, the
estimation of $\psi'$ does not lead to a loss in the rate. As seen in
\eqref{eqDefKe}, we need $\hat\gamma$ to estimate $k_e$.
%
%
%th3.4 #&#
\begin{thmm}\label{thmupperBoundK}
Let $s\geq1, R,\bar\alpha>0$ and assume $Ae^{-A}\lesssim\Delta^2$ as
well as $\Delta(\|\delta_j\|_{l^2}^2+\|(x_j\delta_j)_j\|
_{l^2}^2)\lesssim\|\delta\|_{l^\infty}^2$. Using an estimator $\hat
\gamma$ which satisfies $\sup_{\mathcal P}\mathbb E_{\mathcal
P}[|\hat
\gamma-\gamma|^2]<\infty$ and choosing the cut-off value $U_{\bar
\alpha
}:=\varepsilon^{-2/(2s+2T\bar\alpha+5)}$, we obtain for the risk of
$\hat k_e$ the uniform convergence rate
\[
\sup_{\mathcal P=(\gamma,k)\in\mathcal H_s(R,\bar\alpha)}\mathbb
E_{\mathcal P}\bigl[\|\hat k_e-k_e
\|_{L^2}^2\bigr]^{1/2} \lesssim
\varepsilon^{2s/(2s+2T\bar\alpha+5)}.\vadjust{\goodbreak}
\]
\end{thmm}
%
%
%re3.5 #&#
\begin{remark}
The convergence rates in the Theorems \ref{thmupperBoundsParameters}
and \ref{thmupperBoundK} are minimax optimal up to a logarithmic
factor, which is shown in the supplementary article \cite
{trabs2011Supplement}.
\end{remark}

%To establish asymptotic lower bounds for the convergence rates, we
%consider the continuous white noise model
% \di Z_{\mathcal P}(x)=\mathcal O_{\mathcal P}(x)\di x+\frac{1}{\sqrt
%N}\lambda_N(x) \di W(x), x\in[-A_N,A_N],
%with a two sided Brownian motion $W$, an option function $\mathcal O_{
%apply the results of \cite{brownLow1996} to show the asymptotic
%equivalence of the above considered regression model
%following lower bounds. A detailed description of these results is
%given in the supplemental article \cite{trabs2011Supplement}.
% Let $s\in\mathbb N, s\geq2, R,\bar\alpha>0$ and $j=0,\ldots,s-2,$.
%Then there exists a $\beta\geq0$ such that the following asymptotic
%risk lower bounds hold
% \begin{eqnarray*}
% \inf_{\hat\gamma}\sup_{\mathcal P\in\mathcal G_s(R,\bar
% &\gtrsim\big(\varepsilon(\log\varepsilon^{-1})^{-\beta/2}
% \inf_{\hat\alpha_j}\sup_{\mathcal P\in\mathcal G_s(R,
% &\gtrsim\big(\varepsilon(\log\varepsilon^{-1})^{-\beta/2}
% \inf_{\hat k_e}\sup_{\mathcal P\in\mathcal H_s(R,\bar
% &\gtrsim\big(\varepsilon(\log\varepsilon^{-1})^{-\beta/2}
% \end{eqnarray*}
% where the infimum is taken over all estimators, i.e. all measurable
%functions of the observation $Z$. The bound for $k_e$ holds for $s=1$
%as well.
%The precise $\beta$ depends on the function $\lambda$. It would be
%zero if $\lambda/\sqrt{N}>c$ holds for some $c>0$. In the case $
%convergence rates in the last section. This might be necessary owing
%to the conditions $\Delta\|(\delta_j)\|_{l^2}^2\lesssim\|(\delta_j)
%respectively.

%s4 #&#
\section{Adaptation}\label{secadaptation}
The convergence rate of our estimation procedure depends on the bound
$\bar\alpha$ of the true but unknown $\alpha\in\mathbb R_+$. Therefore,
we construct an $\alpha$-adaptive estimator. For simplicity we
concentrate on the estimation of $\alpha$ itself whereas the results
can be easily extended to $\gamma$, $\alpha_j, j=1,\ldots,s-2$, and
$k_e$. In this section, we will require the following assumption.
%
%
%as4 #&#
\begin{assumption}\label{assparameters}
Let $R>0$, $s\geq2$ and $\alpha\in[0,\bar\alpha]$ for some maximal
$\bar\alpha>0$. Furthermore, we suppose $\mathrm{e}^{-A}\lesssim\Delta^2$ and
$\Delta\| \delta\|_{l^2}^2\lesssim\|\delta\|_{l^\infty}^2$.
\end{assumption}
These conditions only recall the setting in which the convergence rates
of our parameter estimators were proven. Given a consistent
preestimator $\hat\alpha_{\mathrm{pre}}$ of $\alpha$, let $\tilde\alpha_0$ be
the estimator using the data-driven cut-off value and the trimming parameter
%
%
%e4.1 #&#
%e4.2 #&#
\begin{eqnarray}
\tilde U&:=&U_{\hat\alpha_{\mathrm{pre}}}:=\varepsilon^{-2/(2s+2T\hat\alpha
_{\mathrm{pre}}+1)} \quad\mbox{and}
\label{equEpsilon}
\\
\tilde\kappa(u)&:=&\kappa_{\bar\alpha_{\mathrm{pre}}}(u):=\cases{ \displaystyle\tfrac{1}{2}\mathrm{e}^{-TR},&\quad
$|u|<1,$\vspace*{2pt}
\cr
\displaystyle\tfrac{1}{2}C_{\bar\alpha_{\mathrm{pre}}}|u|^{-T\bar\alpha_{\mathrm{pre}}},&\quad
$|u|\geq1,$} \label{eqkappaEpsilon}
\end{eqnarray}
respectively, with $\bar\alpha_{\mathrm{pre}}:=\hat\alpha_{\mathrm{pre}}+|\log
\varepsilon
|^{-1}$. If $\hat\alpha_{\mathrm{pre}}$ is sufficiently concentrated around the
true value, the adaptation does not lead to losses in the rate as the
following proposition shows. Note that the condition $\tilde\alpha_0\in
[0,\bar\alpha]$ is not restrictive since any estimator $\hat\alpha$ of
$\alpha\in[0,\bar\alpha]$ can be improved by using $(0\vee\hat
\alpha
)\wedge\bar\alpha$ instead.
%
%
%pr4.1 #&#
\begin{prop}\label{propadaptiveEstimator}
On Assumption \ref{assparameters} let $\hat\alpha_{\mathrm{pre}}$ be a
consistent estimator which is independent of the data $O_j, j=1,\ldots
,N,$ and fulfills for $\varepsilon\to0$ the inequality
%
%
%e4.3 #&#
\begin{equation}
\label{eqinequalityAlphaPre} \mathbb P\bigl(|\hat\alpha_{\mathrm{pre}}-\alpha|
\geq|\log\varepsilon|^{-1}\bigr)\leq d\varepsilon^2
\end{equation}
with a constant $d\in(0,\infty)$. Furthermore, we suppose $\tilde
\alpha_0\in[0,\bar\alpha]$ almost surely. Then $\tilde\alpha_0$
satisfies the
asymptotic risk bound
\[
\sup_{\mathcal P\in\mathcal G_s(R,\alpha)}\mathbb E_{\mathcal
P,\hat\alpha_{\mathrm{pre}}}\bigl[|\tilde\alpha_0-
\alpha|^2\bigr]^{1/2} \lesssim\varepsilon^{2(s-1)/(2s+2T\alpha+1)},
\]
where the expectation is taken with respect to the common distribution
$\mathbb P_{\mathcal P,\hat\alpha_{\mathrm{pre}}}$ of the observations
$O_1,\ldots
,O_N$ and the preestimator $\hat\alpha_{\mathrm{pre}}$.
\end{prop}
To use $\hat\alpha_0$ on an independent sample as preestimator, we
establish a concentration result for the proposed procedure. We require
$(\varepsilon_j)$ to be uniformly sub-Gaussian (see, e.g., van~de Geer \cite
{vandegeer2000}). That means there\vadjust{\goodbreak} are constants $C_1,C_2\in(0,\infty
)$ such that the following concentration inequality holds for all $t,
N>0$ and $a_1,\ldots, a_N\in\mathbb R$
%
%
%e4.4 #&#
\begin{equation}
\label{eqconcentrationEpsilon} \mathbb P \Biggl( \Biggl|\sum
_{j=1}^Na_j\varepsilon_j \Biggr|
\geq t \Biggr)\leq C_1\exp\biggl(-C_2\frac{t^2}{\sum_{j=1}^Na_j^2}
\biggr).
\end{equation}

%
%pr4.2 #&#
\begin{prop}\label{propconcentrationAlpha}
Additionally to Assumption \ref{assparameters} let $(\varepsilon_j)$
be uniformly sub-Gaussian fulfilling~\eqref{eqconcentrationEpsilon}.
Then there is a constant $c>0$ and for all $\kappa>0$ there is an
$\varepsilon_0 \sim\kappa^{(2s+2T\bar\alpha+1)/(2s-2)}$, such that for
all $\varepsilon<\varepsilon_0\wedge1$ the estimator $\hat\alpha_0$ satisfies
%
%
%e4.5 #&#
\begin{equation}
\label{eqexponentialConcentrationAlpha} \mathbb P\bigl(|\hat\alpha_0-
\alpha|\geq\kappa\bigr)\leq\bigl((7N+1)C_1+2\bigr)\exp\bigl(-c\bigl(
\kappa^2\wedge\kappa^{1/2}\bigr)\varepsilon^{-(s-1)/(2s+2T\bar
\alpha
+1)}
\bigr).
\end{equation}
\end{prop}
Concentration \eqref{eqexponentialConcentrationAlpha} is stronger than
needed in Proposition \ref{propadaptiveEstimator}. To apply the
proposed estimation procedure, let $S_{\mathrm{pre}}$ and $S$ be two independent
samples with noise levels $\varepsilon_{\mathrm{pre}}$ and $\varepsilon$ as well
as sample sizes $N_{\mathrm{pre}}$ and $N$, respectively. Using $S_{\mathrm{pre}}$ for
the estimator $\hat\alpha_{\mathrm{pre}}$, we construct adaptively $\tilde
\alpha_0$ on $S$. We suppose $N_{\mathrm{pre}}$ grows at most polynomial in
$\varepsilon_{\mathrm{pre}}$, that is $N_{\mathrm{pre}}\lesssim\varepsilon_{\mathrm{pre}}^{-p}$
holds for some $p>0$, cf. \cite{trabs2011Supplement}. To satisfy
\eqref{eqinequalityAlphaPre}, it is sufficient if there exists a power
$q>0$, which can be arbitrary small, such that $\varepsilon_{\mathrm{pre}}\sim
\varepsilon^q$ owing to the exponential inequality \eqref
{eqexponentialConcentrationAlpha}. Using $\varepsilon^2\gtrsim
A_N/N\geq1/N$, we estimate
\[
\frac{N_{\mathrm{pre}}}{N}\lesssim\varepsilon_{\mathrm{pre}}^{-p}
\varepsilon^{2}\sim\varepsilon^{2-pq}\to0
\]
for $q<2/p$. Thus, relatively to all available data the necessary
number of observations for the preestimator tends to zero.

%s5 #&#
\section{Discussion and application}\label{sec5}

%s5.1 #&#
\subsection{Numerical example}\label{subsecimplementation}
We apply the proposed estimation procedure to the variance gamma model.
In view of the empirical study \cite{madan1998} we choose the
parameters $\nu\in\{0.05,0.1,0.2,0.5\}, \sigma=1.2$ and $\theta=-0.15$.
the martingale condition \eqref{eqmartingalCondition} yields then
$\gamma=\frac{1}{\nu}\log(1-\theta\nu-\sigma^2\nu/2)$.
According to the different choices of $\nu$, we set $\bar\alpha=40$ as
maximal value of $\alpha$.

The deterministic design of the sample $\{x_1,\ldots,x_N\}$ is
distributed normally with mean zero and variance $1/3$. The
observations $O_j$ are computed from the characteristic function
$\varphi_T$ using the fast Fourier transform method \cite{carrMadan1999}. The additive noise consists of normal centered random
variables with variance $|\delta\mathcal O(x_j)|^2$ for some $\delta
>0$.

We estimate $q\in\{\gamma,\alpha_0,\alpha_1, k\}$. Hence, we need
$s\geq
4$, see Corollary 8 in \cite{trabs2011Supplement}. By self-decomposablity of
the model we apply the rearranged estimator $\hat k^*$ given by \eqref
{eqRearranged}. We use maturity $T=0.25$, interest $r=0.06$, smoothness
$s=6$, sample size $N=100$ and noise level $\delta=0.01$, which
generates values of $\varepsilon$ on average $0.168$. The results of
1000 Monte Carlo simulations are summarized in Tables \ref
{tabadaptation} and \ref{tabn100delta001}.
%
%
%t1 #&#
\begin{table}
\tablewidth=270pt
\caption{Risk of estimating $\alpha$ with oracle (\textit{middle
column}) and adaptive (\textit{right column}) choice of cut-off value
$U$ in simulated variance gamma model with $\nu\in\{0.05,0.1,0.2,0.5\},
\sigma=1.2, \theta=-0.15$}\label{tabadaptation}
\begin{tabular*}{270pt}{@{\extracolsep{\fill}}lll@{}}
\hline
$\alpha$ & \multicolumn{1}{l}{$\mathbb E[|\hat\alpha_0-\alpha
|^2]^{1/2}$} & \multicolumn{1}{l@{}}{$\mathbb E[|\tilde\alpha_0-\alpha
|^2]^{1/2}$} \\
\hline
40 & 20.7998 & 23.3589 \\
20 & \phantom{0}5.8362 & \phantom{0}7.7724 \\
10 & \phantom{0}1.0505 & \phantom{0}2.4534 \\
\phantom{0}4 & \phantom{0}0.1729 & \phantom{0}1.1158 \\
\hline
\end{tabular*} \vspace*{-3pt}
\end{table}
%

%
%t2 #&#
\begin{table}
\caption{Risk of estimating the parameters $\gamma, \alpha, \alpha_1$
and the $k$-function with oracle (\textit{middle column}) and adaptive
(\textit{right column}) choice of the cut-off value $U$ in simulated
variance gamma model ($\nu=0.2, \sigma=1.2, \theta=-0.15$)}\label{tabn100delta001}
\begin{tabular*}{\textwidth}{@{\extracolsep{\fill}}llll@{}}
\hline
& \multicolumn{1}{l}{$q$} & \multicolumn{1}{l}{$\mathbb E[|\hat
q-q|^2]^{1/2}$} & \multicolumn{1}{l@{}}{$\mathbb E[|\tilde
q-q|^2]^{1/2}$}\\
\hline
$\gamma$ & \phantom{0$-$}0.1408 & \phantom{0}0.0065 & \phantom{0}0.0126 \\
$\alpha_0$ & \phantom{$-$}10.0000 & \phantom{0}1.0505 & \phantom{0}2.4534 \\
$\alpha_1$ & $-94.1667$ & 32.1016 & 77.5311 \\
\hline
& \multicolumn{1}{l}{$\|q\|_{L^2}^{1/2}$} & \multicolumn
{1}{l}{$\mathbb E[\|\hat q-q\|_{L^2}^2]^{1/2}$} & \multicolumn
{1}{l@{}}{$\mathbb E[\|\tilde q-q\|_{L^2}^2]^{1/2}$} \\
\hline
$k_e$ & 0.9556 & 0.4075 & 0.5602 \\
\hline
\end{tabular*}  \vspace*{-3pt}
\end{table}
In order to apply the estimation procedure, we need to choose the
tuning parameters. Owing to the typically unknown smoothness $s$, let
the weight functions satisfy Assumption \ref{assweightFunctions} for
some large value $s_{\mathrm{max}}$. The weights for the parameters can be
chosen as polynomial whereas $W_k$ is taken as a polynomial times a
smooth function with support $[-1,0]$. The trimming parameter~$\kappa$
is included mainly for theoretical reasons and is not important to the
implementation. The most crucial point is the choice of the cut-off
value $U$. For $\hat q$ we implement the oracle method
$U=\argmin_{V\geq
0} |\hat q(V)-q|$ and an adaptive estimator $\tilde q$ based on the
construction of Section \ref{secadaptation} with sample size
$N_{\mathrm{pre}}=25$ for $\hat\alpha_{\mathrm{pre}}$.\looseness=-1

%s5.2 #&#
\subsection{Discussion}
Due to the nonparametric setting, our estimators converge more slowly
than with $\sqrt{n}$ rate as in parametric models \cite
{eberlein2012,eberlein1998,madan1998}. Although the studied
estimation problem is only mildly ill-posed compared with classical
nonparametric regression models and thus the polynomial rates are
faster than in nonparametric models with $\sigma>0$ which achieve
logarithmic rates only \cite{reiss12006}. In order to understand the
convergence rate of the estimators for $\gamma$ and $\alpha_j$ better,
we rewrite equation \eqref{eqrepPsiPrim} in the distributional sense,
denoting the Dirac distribution at zero by $\delta_0$, and
differentiate representation \eqref{eqrepresentationPsi}
\[
\psi'(u)=\mathcal F (\mathrm{i}\gamma\delta_0+\mathrm{i}k_e
) (u)=\mathrm{i}\gamma-\sum_{j=0}^{s-2}\mathrm{i}^jj!
\alpha_ju^{-j-1}+\rho'(u),\qquad  u\in\mathbb R
\setminus\{0\}.\vadjust{\goodbreak}
\]
Hence, $\psi'$ can be seen as Fourier transform of an s-times weakly
differentiable function and estimating $\gamma$ from noisy observations
of $\psi'$ corresponds to a nonparametric regression with regularity s.
Since dividing by $u$ on the right-hand side of the above equation
corresponds to taking the derivative in the spatial domain, the
estimation of $\alpha_j$ is similar to the estimation of the $(j+1)$th
derivative in a regression model. The convergence rate of $k_e$ is in
line with the results of Belomestny and Rei{\ss} \cite{reiss12006} for $\sigma=0$ since their
rate equals ours in the compound Poisson case $\alpha=0$.

For $\hat k_e$, the degree of ill-posedness is given by $T\alpha+2$.
This can be seen analytically by observing that the noise is governed
by $u^2|\varphi_T(u-\mathrm{i})|^{-1}$, which grows with rate $T\alpha+2$. From
a statistical point of view a higher value of $\alpha$ leads to a more
active L\'evy process and hence, it is harder to distinguish the small
jumps of the process from the additive noise. The influence of the time
to maturity $T$ on the convergence rates is an interesting deviation
from the analysis of~Belomestny and Rei{\ss} \cite{reiss12006}. The simulation shown in Table
\ref{tabadaptation} demonstrates the improvement of the estimation for
small the values of $\alpha$. The estimators $\hat\gamma$ and $\hat
k_e$ provide a complete calibration of the model. Although, estimating
the $k$-function at zero is most important and thus additional
information through $\hat\alpha_j$ are crucial. Table \ref
{tabn100delta001} contains simulation results for the estimators $\hat
q$ and $\tilde q$, $q\in\{\gamma, \alpha_0, \alpha_1, k_e\}$,
corresponding to oracle and $\alpha$-adaptive cut-off values,
respectively. This adaptation to $\alpha$ is a first step to a
data-driven procedure and should be developed further.

Since the estimating equation \eqref{eqrepPsi2} holds for all L\'evy
processes with finite variation, the proposed estimator $\hat k_e$ can
be more generally understood as estimator of $xe^x\nu(\mathrm
{d}x)$. Thus,
the estimation procedure can be applied to exponential L\'evy models
with Blumenthal--Getoor index larger than zero, for example, tempered
stable processes. However, the analysis of the convergence rates does
not carry over to more general L\'evy processes since the polynomial
decay of the $\varphi_T$, which is guaranteed by property \eqref
{propK}, is essential for our proofs. Moreover, if $k$ has no bounded
variation the behavior of the L\'evy density at zero needs different
methods and should be studied further. For instance, Belomestny \cite{belomestny2010}
discusses the estimation of the fractional order for
regular L\'evy models of exponential type.\looseness=-1

Even if the practitioner prefers specific parametric models that might
achieve smaller errors and faster rates, the nonparametric method
should be used as a goodness-of-fit test against model
misspecification. To construct such tests, confidence sets need to be
studied which is done by S{\"o}hl~\cite{soehl2012} in the framework of L\'evy
processes with finite activity. Based on this asymptotic analysis,
S{\"o}hl and Trabs \cite{soehlTrabs2012} construct confidence intervals in the
self-decomposable model.

%t3 #&#
\begin{table}
\tablewidth=200pt
\caption{Adaptive estimation based on DAX options from 29 May 2008 with
time to maturity $T$ and $N+N_{\mathrm{pre}}$ observations}\label{tabrealData}
\begin{tabular*}{200pt}{@{\extracolsep{\fill}}ld{3.3}d{3.3}@{}}
\hline
$T$ & 0.314 & 0.567 \\
$r$ & 0.045 & 0.044 \\
$N_{\mathrm{pre}}$ & 20 & 21 \\
$N$ & 81 & 85 \\[3pt]
$\tilde\gamma$ & 0.101 & 0.344 \\
$\tilde\alpha_0$ & 34.848 & 23.600 \\
$\tilde\alpha_1$ & 239.348 & 147.699 \\
\hline
\end{tabular*}
\end{table}

\vspace*{-3pt}
%s5.3 #&#
\subsection{Real data example}
\vspace*{-3pt}%

\begin{figure}[b]

\includegraphics{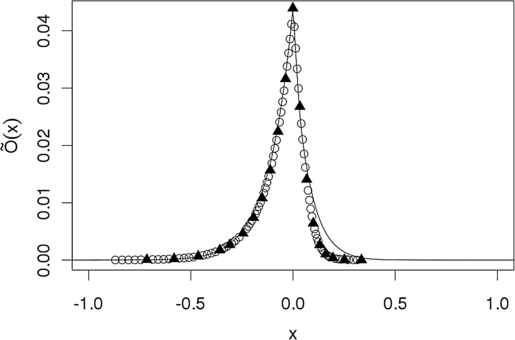}

\caption{Observed DAX option prices (\textit{points and triangles})
from 29 May 2008 with fixed maturity $T=0.314$ and different log strike
prices $x$ as well as the option function generated from the estimated
model (\textit{solid line}).}\label{figoptionFunction3months}
\end{figure}

\begin{figure}

\includegraphics{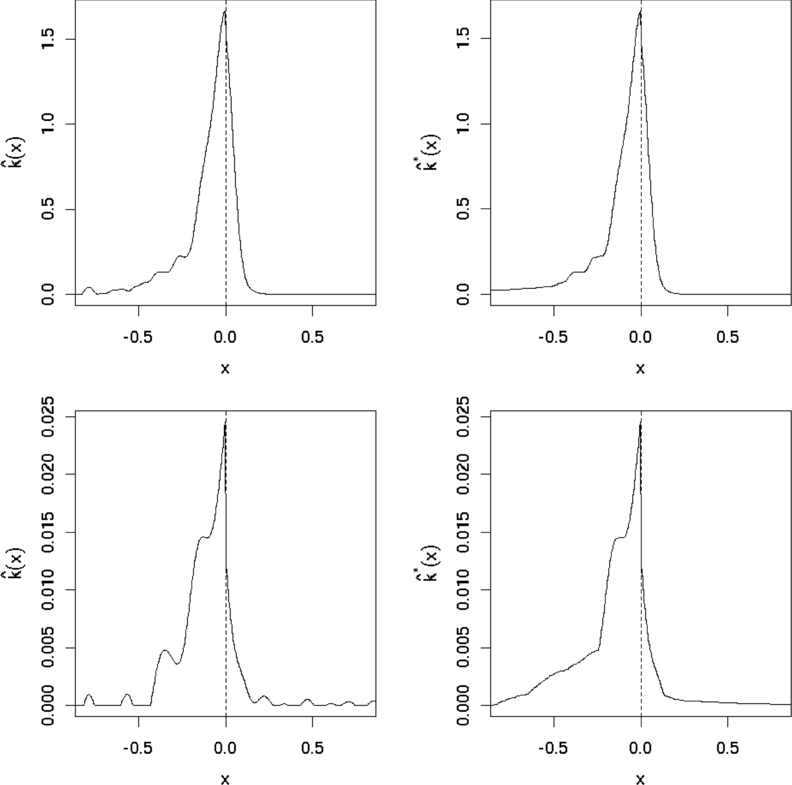}

\caption{Using ODAX data from 29 May 2008 with three (\textit{top}) and
six (\textit{bottom}) months to maturity, estimation of the function
$k$ with (\textit{right}) and without (\textit{left})
rearrangement.}\label{figkFunction3months}
\end{figure}

We apply our estimation method to a data set from the Deutsche B\"orse
database Eurex.\footnote{Provided through the Collaborative Research
Center 649 ``Economic Risk''.} It consists of settlement prices\vadjust{\goodbreak} of put
and call options on the DAX index with three and six months to maturity
from 29 May 2008. The sample sizes are 101 and 106, respectively. The
interest rate is chosen according to the put-call parity. The
sub-sample for the preestimator consists of every fifth strike while
the main estimation is done from the remaining data points. By a rule
of thumb, the bid-ask spread is chosen as 1\% of the option prices.
Therefore, we get noise levels $\varepsilon$ with values 0.0138 and
0.069 for the two maturities, respectively. Table \ref{tabrealData}
shows the result of the proposed method. As one would expect, the jump
activity is smaller for a longer time to maturity. The estimator $\hat
k(x)=\mathrm{e}^{-x}\hat k_e(x)$ as well as the rearranged estimator $\hat k^*$
are presented in Figure~\ref{figkFunction3months}. In Figure \ref
{figoptionFunction3months}, the calibrated model is used to generate
the option function in the case of three months to maturity, where the
data points used for the preestimator are marked with triangles in the
figure. For a comparison of the outcome of our estimation procedure
with the spectral calibration of Belomestny and Rei{\ss} \cite{reiss12006}, we refer to
S{\"o}hl and Trabs \cite{soehlTrabs2012}.

%f1 #&#

%
%
%f2 #&#

%s6 #&#
\section{Proofs}\label{sec6}

%s6.1 #&#
\subsection{Proof of the upper bounds}\label{subsecproofUpperBoundParameters}
Let us recall some results of \cite{reiss12006}: Because of the
$B$-spline interpolation we obtain $\mathcal O_l(x):=\mathbb E[\tilde
{\mathcal O}(x)]=\sum_{j=1}^{N}\mathcal O(x_j)b_j(x)+\beta_0(x),x\in
\mathbb R$. Furthermore, the decomposition of the stochastic error
$\tilde\psi-\psi$ in a linearization $\mathcal L$ and a remainder
$\mathcal R$,
\begin{eqnarray*}
\mathcal L(u):=T^{-1}\varphi_T(u-\mathrm{i})^{-1}(\mathrm{i}-u)u
\mathcal{F}(\tilde\mathcal{O}-\mathcal{O})(u), \qquad \mathcal R(u):=\tilde\psi(u)-\psi
(u)-\mathcal L(u),
\end{eqnarray*}
$u\in\mathbb R,$ has the following properties.
%
%
%pr6.1 #&#
\begin{prop}\label{propresultsReiss}
\textup{(i)} Under the hypothesis $\mathrm{e}^{-A}\lesssim\Delta^2$ we obtain $\sup
_{u\in\mathbb R}|\mathbb E[\mathcal{F}\tilde\mathcal{O}(u)-\mathcal
{FO}(u)]|=\sup_{u\in\mathbb R}|\mathcal{FO}_l(u)-\mathcal
{FO}(u)|\lesssim\Delta^2$
uniformly over all L\'evy triplets satisfying Assumption \ref{assmoment}.\vspace*{-6pt}
\begin{longlist}[(ii)]
\item[(ii)] If the function $\kappa\dvt \mathbb R\to\mathbb R_+$ satisfies
\eqref
{eqconditionKappa}, then for all $u\in\mathbb R$ the remainder is
bounded by $|\mathcal R(u)|\leq T^{-1}\kappa(u)^{-2}(u^4+u^2)|\mathcal
F(\tilde\mathcal{O}-\mathcal{O})(u)|^2$.
\end{longlist}
\end{prop}
\subsubsection*{Upper bound for $\gamma$ and $\alpha_j$ (Theorem
\protect\ref{thmupperBoundsParameters})}
Since Theorem \ref{thmupperBoundsParameters} can be proven
analogously to Theorem 4.2 in \cite{reiss12006}, we only sketch the
main steps. Note that in $\mathcal G_s(R,\bar\alpha)$ we can bound
uniformly $|u^{s-1}\rho(u)|$ in representation \eqref
{eqrepresentationPsi}, cf. Lemma 9 in \cite{trabs2011Supplement}. Let us
consider $\gamma$ first. The definition of $\hat\gamma$ and
$w_\gamma^U$, the decomposition of $\tilde\psi$ and representation
\eqref
{eqrepresentationPsi} yield
\begin{eqnarray*}
\hat\gamma=\int_{-U}^U\im\bigl(\tilde\psi(u)
\bigr)w_\gamma^U(u)\,\mathrm{d}u =\gamma+\int
_{-U}^U\im\bigl(\rho(u)+\mathcal L(u)+\mathcal
R(u)\bigr)w_\gamma^U(u)\,\mathrm{d}u.
\end{eqnarray*}
Hence, we obtain
\begin{eqnarray*}
\mathbb E\bigl[|\hat\gamma-\gamma|^2\bigr]&\leq&3 \biggl|\int
_{-U}^U\rho(u)w_\gamma^U(u)
\,\mathrm{d}u\biggr |^2 +3\mathbb E \biggl[ \biggl|\int_{-U}^U
\mathcal L(u)w_\gamma^U(u)\,\mathrm{d}u \biggr|^2
\biggr]
\\
&&{} +3\mathbb E \biggl[ \biggl|\int_{-U}^U\mathcal
R(u)w_\gamma^U(u)\, \mathrm{d} u \biggr|^2 \biggr],
\end{eqnarray*}
where all three summands can be estimated separately. The first one is
a deterministic error term. It can be estimated using the decay of
$\rho
(u)$ and the weight function property~\eqref{eqweightfunctionProperty}:
\[
\biggl|\int_{-U}^U\rho(u)w_\gamma^U(u)
\,\mathrm{d}u \biggr|\lesssim\int_{-U}^UU^{-(s+1)}\bigl|
\rho(u)u^{s-1}\bigr|\,\mathrm{d}u\lesssim U^{-s}.
\]
A bias-variance decomposition, with the definition $\var(Z):=\mathbb
E[|Z-\mathbb E[Z]|^2]$, of the linear error term yields
\begin{eqnarray*}
&&\mathbb E \biggl[ \biggl|\int_{-U}^U\mathcal
L(u)w_\gamma^U(u)\,\mathrm{d}u \biggr|^2 \biggr]\\
&&\quad=
\biggl|\int_{-U}^U\frac{(\mathrm{i}-u)u}{T\varphi_T(u-\mathrm{i})}\mathbb E\bigl[
\mathcal{F}(\tilde\mathcal{O}-\mathcal{O})(u)\bigr]w_\gamma^U(u)\,\mathrm{d}u
\biggr|^2
\\
&&\qquad{}+\var\biggl(\int_{-U}^U\frac{(\mathrm{i}-u)u}{T\varphi_T(u-\mathrm{i})}
\mathcal{F}\tilde\mathcal{O}(u)w_\gamma^U(u)\,\mathrm{d}u \biggr) =:
\mathcal L_b^2+\mathcal L_v.
\end{eqnarray*}
Using the approximation result in Proposition \ref{propresultsReiss},
the bound of $|\varphi_T(u-\mathrm{i})|^{-1}$ given by $\kappa^{-1}$ and
property \eqref{eqweightfunctionProperty}, we infer the estimate of
the bias term
\begin{eqnarray*}
|\mathcal L_b| \lesssim\Delta^2U^{-(s+1)}\int
_{-U}^U\bigl|\varphi_T(u-\mathrm{i})\bigr|^{-1}|u|^{s+1}
\,\mathrm{d}u\lesssim\Delta^2U^{T\bar\alpha+1}.
\end{eqnarray*}
%
%For the variance part we make use of the properties of the linear
%spline functions $b_k$: For $k=1,\ldots,N$ they satisfy %$
%$|x_{k+1}-x_{k-1}|\leq2\Delta$ follows
% \|\mathcal Fb_k\|_{L^2}=\sqrt{2\pi}\|b_k\|_{L^2}\leq\sqrt{4\pi\Delta}
% \mbox{and} \|\mathcal Fb_k\|_\infty\leq\|b_k\|%_{L^1}\leq2
For the variance part, we make use of the properties of the linear
spline functions $b_k$ as well as $\supp(w_\gamma^U)\subseteq[-U,U]$
and the independence of $(\varepsilon_k)$. We estimate ($\cov
(Y,Z):=\mathbb E[(Y-\mathbb E[Y])\overline{(Z-\mathbb E[Z])}]$)
\begin{eqnarray*}
\mathcal L_v&=&\int_{-U}^U\int
_{-U}^U\cov\biggl(\frac
{(\mathrm{i}-u)u}{T\varphi_T(u-\mathrm{i})}\mathcal{F}
\tilde\mathcal{O}(u),\frac{(\mathrm{i}-v)v}{T\varphi_T(v-\mathrm{i})}\mathcal{F}\tilde\mathcal{O}(v)
\biggr)w_\gamma^U(u)w_\gamma^U(v)
\,\mathrm{d}u\,\mathrm{d}v
\\
&=&\sum_{k=1}^N\delta_k^2
\biggl|\int_{-U}^U\frac{(\mathrm{i}-u)u}{T\varphi_T(u-\mathrm{i})}\mathcal
Fb_k(u)w_\gamma^U(u)\,\mathrm{d}u
\biggr|^2 \lesssim\Delta\|\delta\|_{l^\infty}^2U^{2T\bar\alpha+1}.
\end{eqnarray*}
To estimate the remaining term $\mathcal R$, we use Proposition \ref
{propresultsReiss}, the property \eqref{eqweightfunctionProperty} of
$w_\gamma^U$ and the choice of $\kappa$. In addition the independence
of $(\varepsilon_k)$ and the uniform bound of their fourth moments
comes into play.
\begin{eqnarray*}
&&\mathbb E \biggl[\biggl |\int_{-U}^U\mathcal
R(u)w_\gamma^U(u)\,\mathrm{d}u \biggr|^2 \biggr]
\\
&&\quad\lesssim\int_{-U}^U\int_{-U}^U
\bigl(\bigl\|\mathcal{F}(\mathcal{O}_l-\mathcal{O})\bigr\|_\infty^4
+\mathbb E \bigl[\bigl |\mathcal{F}(\tilde\mathcal{O}- \mathcal{O}_l) (u)\mathcal{F}(\tilde\mathcal{O}-\mathcal{O}_l) (v)\bigr |^2 \bigr] \bigr)\\
&&\hspace*{56pt}{}\times \frac{u^4w_\gamma^U(u)v^4w_\gamma
^U(v)}{\kappa(u)^2\kappa(v)^2}\,
\mathrm{d}u\,\mathrm{d}v
\\
&&\quad\lesssim\biggl(\Delta^4\int_{-U}^U
\frac{u^4w_\gamma^U(u)}{\kappa
(u)^2}\,\mathrm{d}u \biggr)^2+ \Biggl(\int
_{-U}^U\sum_{k=1}^N
\delta_k^2\bigl|\mathcal Fb_k(u)\bigr|^2
\frac{u^4w_\gamma^U(u)}{\kappa(u)^2}\, \mathrm{d}u \Biggr)^2
\\
&&\quad\lesssim\biggl(\Delta^4U^{-(s+1)}\int_{-U}^U
\kappa(u)^{-2}|u|^{s+3}\,\mathrm{d} u \biggr)^2+
\biggl(\Delta^2\|\delta\|_{l^2}^2U^{-(s+1)}
\int_{-U}^U\kappa(u)^{-2}|u|^{s+3}
\,\mathrm{d}u \biggr)^2
\\
&&\quad\lesssim U^{4T\bar\alpha+6}\bigl(\Delta^8+\Delta^4\|
\delta\|_{l^2}^4\bigr).
\end{eqnarray*}
Therefore, the total risk of $\hat\gamma$ is of order
\[
\mathbb E\bigl[|\hat\gamma-\gamma|^2\bigr]\lesssim
U^{-2s}+U^{2T\bar\alpha
+1}\bigl(\Delta^4U+\Delta\|\delta
\|_{l^\infty}^2\bigr)+U^{4T\bar\alpha
+6}\bigl(\Delta^8+
\Delta^4\|\delta\|_{l^2}^4\bigr)
\]
uniformly over $\mathcal G_s(R,\bar\alpha)$. Since the explicit choice
of $U=U_{\bar\alpha}=\varepsilon^{-2/(2s+2T\bar\alpha+1)}$ fulfills
$U\lesssim\Delta^{-1}$ and $\Delta\|\delta\|_{l^2}^2\lesssim\|
\delta\|_{l^\infty}^2$ holds by assumption, this bound simplifies to
\[
\mathbb E\bigl[|\hat\gamma-\gamma|^2\bigr]\lesssim
U^{-2s}+U^{2T\bar\alpha
+1}\varepsilon^2+U^{4T\bar\alpha+6}
\varepsilon^4.
\]
Here $U_{\bar\alpha}$ balances the trade-off between the first and the
second term whereby the third term is asymptotically negligible. We
obtain the claimed rate.

For $\alpha_j,j=0,\ldots,s-2,$ the only difference to the analysis for
$\hat\gamma$ is the rescaling factor of $w_{\alpha_j}^U$ in \eqref
{eqweightfunctionProperty}. Since its square appears in front of every
term, we verify
\begin{eqnarray*}
\mathbb E\bigl[|\hat\alpha_j-\alpha_j|^2
\bigr]&\lesssim& U^{-2(s-1-j)}+U^{2T\bar
\alpha+2j+3}\bigl(\Delta^4U+
\Delta\|\delta\|_{l^\infty}^2\bigr) +U^{4T\bar\alpha+2j+8}\bigl(
\Delta^8+\Delta^4\|\delta\|_{l^2}^4
\bigr)
\\[-2pt]
&\lesssim& U^{-2(s-1-j)}+U^{2T\bar\alpha+2j+3}\varepsilon^2+U^{4T\bar
\alpha+2j+8}
\varepsilon^4.
\end{eqnarray*}
%
% The explicit choice of $U=U_{\bar\alpha}$ implies the result.

\subsubsection*{Upper bound for $k_e$ (Theorem \protect\ref
{thmupperBoundK})}
Similarly to the uniform bound of the bias of $\mathcal{F}\tilde\mathcal{O}$ in
Proposition \ref{propresultsReiss}, the following lemma holds true. It
can be proved analogously to Proposition 1 in \cite{reiss12006} and thus we
omit the details.
%
%
%le6.2 #&#
\begin{lemma}\label{lemxbBound}
Assuming $Ae^{-A}\lesssim\Delta^2$, we obtain $\sup_{u\in\mathbb
R}|\mathbb E[\mathcal F (x(\tilde\mathcal{O}-\mathcal{O})(x)
)(u)]|= \sup_{u\in\mathbb R}|\mathcal F (x(\mathcal O_l-\mathcal O)(x)
)(u)|\lesssim\Delta^2$ uniformly over all L\'evy triplets satisfying
Assumption \ref{assmoment} and $\mathbb E[|X_Te^{X_T}|]\lesssim1$.
\end{lemma}
%
% We follow the lines of the proof of Proposition 6.1 in
% \begin{eqnarray*}
% &\int_{x_1}^{x_N}|x(\mathcal O_l-\mathcal O)(x)\di x|
% \leq\sum_{j=2}^N(|x_{j-1}|\vee|x_j|)\int_{x_{j-1}}^{x_j}|\mathcal
%O_l(x)-\mathcal O(x)|\di x\\
% \leq&\sum_{j\in\{2,\ldots,N\}\setminus\{j_0\}}\int_{x_{j-1}}^{x_j}
%O''(z)|\di z\di y\di x+C_0(|x_{j_0-1}|\vee|x_{j_0}|)\Delta^2\\
% \leq&\|x\mathcal O''(x)\|_{L^1}\Delta^2+\|\mathcal O''\|_{L^1}
% \end{eqnarray*}
% Since the extrapolation errors can be bounded by $4C_2\Delta(A+
% \begin{eqnarray*}
% &\int_{-\infty}^\infty|\mathbb E[x(\tilde\mathcal{O}-\mathcal O)(x)]|
% \leq& 2C_2(Ae^{-A} +e^{-A}) +\|x\mathcal O''(x)\|_{L^1}\Delta^2 +(\|
% \end{eqnarray*}
% It remains to bound $\|x\mathcal O''(x)\|_{L^1}$. Recall from
%P(X_T<x)+f_T(x)-\mathbf1_{\{x>0\}}\big), x\in\mathbb R$.
%Integration by parts yields
% \begin{eqnarray*}
% \int_0^\infty|x\mathcal O''(x)|\di x
% &=\int_0^\infty xe^x\big|\mathbb P(X_T<x)+f_T(x)-1\big|\di x\\
% &\leq\int_0^\infty xe^x\big(1-\mathbb P(X_T<x)\big)\di x+\mathbb E
% &=1-\mathbb P(X_T<0)+\int_0^\infty(x-1)e^xf_T(x)\di x+\mathbb E
% &=\mathbb P(X_T\geq0)+\mathbb E\big[(2|X_T|+1)e^{X_T}\mathbf1_{
% \end{eqnarray*}
% We conclude analogously
% \begin{eqnarray*}
% \int_{-\infty}^0|x\mathcal O''(x)|\di x
% &=\mathbb P(X_T<0)+\mathbb E\big[(2|X_T|+1)e^{X_T}\mathbf1_{\{X_T<0
% \end{eqnarray*}
% Therefore, it holds $\|x\mathcal O''(x)\|_{L^1}\leq2+2\mathbb
%E[|X_Te^{X_T}|]$, which is bounded by assumption.
For convenience, we write $m:=2s-1$ and $w_k:=\mathcal F W_k$ such that
$w_k(u/U)=U\times \mathcal F(W_k(Ux))(u)$. Using $\|f\|_{L^2}^2=\int_{\R
_+}|f(x)|^2\,\mathrm{d}x+\int_{\R_-}|f(x)|^2\,\mathrm{d}x=:\|f\|_{L^2(\R
_+)}^2+\|f\|_{L^2(\R_-)}^2$ for $f\in L^2(\R)$, it is sufficient to
consider the
loss of $\hat k_e$ on $\R_+$. On $\R_-$ one can proceed analogously. We
split the risk into a deterministic error, an error caused by $\hat
\gamma$ and a stochastic error,
\begin{eqnarray*}
&&\mathbb E_{\mathcal P}\bigl[\|\hat k_e-k_e
\|_{L^2(\R_+)}^2\bigr]
\\[-2pt]
&&\quad=\mathbb E_{\mathcal P}\biggl[\biggl\|\mathcal F^{-1} \biggl(\bigl(-\hat
\gamma-\mathrm{i}\tilde\psi'(u)\bigr)w_k\biggl(
\frac{u}{U}\biggr) \biggr)-k_e\biggr\|_{L^2(\R_+)}^2
\biggr]
\\[-2pt]
&&\quad\leq\mathbb E_{\mathcal P} \biggl[\int_{\R_+}3 \biggl|\mathcal
F^{-1} \biggl(\bigl(-\gamma-\mathrm{i}\psi'(u)
\bigr)w_k\biggl(\frac{u}{U}\biggr) \biggr)
(x)-k_e(x) \biggr|^2
\\[-2pt]
&&\hspace*{36pt}\qquad{} +3 \biggl|\mathcal F^{-1} \biggl((\gamma-\hat\gamma) w_k
\biggl(\frac
{u}{U}\biggr) \biggr) (x) \biggr|^2\\[-2pt]
&&\hspace*{36pt}\qquad{}+3 \biggl|\mathcal
F^{-1} \biggl(\bigl(-\mathrm{i}\tilde\psi'(u)+\mathrm{i}
\psi'(u)\bigr)w_k\biggl(\frac
{u}{U}\biggr)
\biggr) (x) \biggr|^2\,\mathrm{d}x \biggr]
\\[-2pt]
&&\quad\le3\int_{\R_+} \biggl|\mathcal F^{-1} \biggl(\mathcal
Fk_e(u)w_k\biggl(\frac
{u}{U}\biggr) \biggr)
(x)- k_e(x) \biggr|^2\,\mathrm{d}x\\[-2pt]
&&\qquad{} +3\mathbb E\bigl[|\hat
\gamma-\gamma|^2\bigr]\int_{\R_+}\bigl|UW_k(Ux)\bigr|^2
\,\mathrm{d}x
\\[-2pt]
&&\qquad{} +3\mathbb E \biggl[\int_{\R_+} \biggl|\mathcal F^{-1}
\biggl( \bigl(\tilde\psi'(u)-\psi'(u)
\bigr)w_k\biggl(\frac{u}{U}\biggr) \biggr) (x) \biggr|^2
\,\mathrm{d}x \biggr]
\\[-2pt]
&&\quad=:D+G+S.
\end{eqnarray*}
The support of $W_k$ yields $G=0$. The deterministic term $D$ can be
estimated in the spatial domain, where we use the local smoothness of
$k_e$. For pointwise convergence rates, this was done in \cite
{belomestny2011}. We decompose using $\supp W_k\subset(-\infty,0]$
\begin{eqnarray*}
D&=&3\int_{\R_+} \bigl|U \bigl(k_e\ast
W_k(U\bullet) \bigr) (x)-k_e(x) \bigr|^2\,
\mathrm{d}x
\\
&\leq&6\int_{\R_+} \biggl|\int_{U}^{\infty}
\bigl(k_e(x+y/U)-k_e(x) \bigr)W_k(-y)\,
\mathrm{d}y \biggr|^2\,\mathrm{d}x
\\
&&{} +6\int_{\R_+} \biggl|\int_0^U
\bigl(k_e(x+y/U)-k_e(x) \bigr)W_k(-y)\,
\mathrm{d} y \biggr|^2\,\mathrm{d}x \\
&=:&6(D_1+D_2).
\end{eqnarray*}
Cauchy--Schwarz's inequality, the estimate $\int_U^\infty|W_k(-y)|\,
\mathrm{d}
y\leq U^{-m}\int_{\mathbb R}|y^{m}W_k(y)|\,\mathrm{d}y\lesssim
U^{-m}$ and
Fubini's theorem yield
\begin{eqnarray*}
D_1&\leq&\int_{\R_+}\int_U^\infty\bigl|W_k(-y)\bigr|
\,\mathrm{d}y\int_U^\infty\bigl|k_e(x+y/U)-k_e(x)
\bigr|^2\bigl|W_k(-y)\bigr|\,\mathrm{d}y\,\mathrm{d}x
\\
&\lesssim& U^{-m}\int_{\R_+}\int
_U^\infty\bigl(\bigl|k_e(x+y/U)\bigr|^2+\bigl|k_e(x)\bigr|^2
\bigr)\bigl|W_k(-y)\bigr|\,\mathrm{d}y\,\mathrm{d}x
\\
&\lesssim& U^{-m}\int_U^\infty\bigl|W_k(-y)\bigr|
\int_{\R
_+}\bigl|k_e(x+y/U)\bigr|^2+\bigl|k_e(x)\bigr|^2
\,\mathrm{d}x\,\mathrm{d}y \lesssim U^{-2m}\|k_e
\|_{L^2}^2.
\end{eqnarray*}
Using a Taylor expansion, we split $D_2$ in a polynomial part and a remainder:
\begin{eqnarray*}
D_2&\leq&2\int_{\R_+} \Biggl|\int_0^U
\Biggl(\sum_{j=0}^{s-1}\frac
{k_e^{(j)}(x)}{j!U^j}y^j
\Biggr)W_k(-y)\,\mathrm{d}y \Biggr|^2\,\mathrm{d}x
\\
&&{} +2\int_{\R_+} \biggl|\int_0^U
\int_x^{x+{y}/{U}}\frac
{k_e^{(s)}(z)(x+y/U-z)^{s-1}}{(s-1)!}\,\mathrm{d}z
W_k(-y)\,\mathrm{d}y \biggr|^2\,\mathrm{d}x
=:2D_{2P}+2D_{2R}.
\end{eqnarray*}
We estimate by $\int_0^Uy^jW_k(-y)\,\mathrm{d}y=-\int_U^\infty
y^jW_k(-y)\,\mathrm{d}y$
for $j=0,\ldots,s-1$
\begin{eqnarray*}
D_{2P}\leq sU^{-2m}\sum_{j=0}^{s-1}
\frac{1}{(j!)^2}\int_{\R
_+}\bigl|k_e^{(j)}(x)\bigr|^2
\,\mathrm{d}x \biggl(\int_{\R_+}\bigl|y^{m}W_k(-y)\bigr|
\, \mathrm{d}y \biggr)^2 \lesssim U^{-2m}\sum
_{j=0}^{s-1}\bigl\|k_e^{(j)}
\bigr\|_{L^2}^2.
\end{eqnarray*}
With twofold usage of Cauchy--Schwarz and with Fubini's theorem we obtain
\begin{eqnarray*}
D_{2R} &=&\int_{\R_+} \biggl|\int_0^U
\int_0^{y/U}\frac
{k_e^{(s)}(x+z)(y/U-z)^{s-1}}{(s-1)!}\,\mathrm{d}z
W_k(-y)\,\mathrm{d}y \biggr|^2\,\mathrm{d}x
\\
&\leq&\int_{\R_+} \biggl(\int_0^U
\biggl(\int_0^{
{y}/{U}}\bigl|k_e^{(s)}(x+z)\bigr|^2
\,\mathrm{d}z \biggr)^{1/2} \\
&&\hspace*{41pt}{}\times\biggl(\int_0^{
{y}/{U}}
\frac
{({y}/{U}-z)^{2s-2}}{((s-1)!)^2}\,\mathrm{d}z \biggr)^{1/2}
\bigl|W_k(-y)\bigr|\,
\mathrm{d}y \biggr)^2\,\mathrm{d}x
\\
&\leq&\int_{\R_+}\int_0^U
\int_0^{{y}/{U}}\bigl|k_e^{(s)}(x+z)\bigr|^2
\, \mathrm{d} z\bigl|W_k(-y)\bigr|\,\mathrm{d}y\int_0^U
\frac
{(y/U)^{2s-1}}{(2s-1)((s-1)!)^2}\bigl|W_k(-y)\bigr|\,\mathrm{d}y\,\mathrm{d}x
\\
&\lesssim& U^{-(2s-1)}\int_0^U\int
_0^{y/U}\int_{\R
_+}\bigl|k_e^{(s)}(x+z)\bigr|^2
\,\mathrm{d}x\,\mathrm{d}z\bigl|W_k(-y)\bigr|\,\mathrm{d}y
\\
&\leq& U^{-(2s-1)}\bigl\|k_e^{(s)}\bigr\|_{L^2}^2
\int_0^U\frac{y}{U}\bigl|W_k(-y)\bigr|
\, \mathrm{d}y \lesssim U^{-2s}.
\end{eqnarray*}
Therefore, we have $D+G\lesssim U^{-2s}+U\mathbb E[|\hat\gamma-\gamma
|^2]$.

To estimate the stochastic error $S$, we bound the term $|\tilde\psi
'(u)-\psi'(u)|$. Let us introduce the notation
\begin{eqnarray*}
\tilde\varphi_T(u-\mathrm{i})&:=&v_{\kappa(u)} \bigl(1+
\bigl(\mathrm{i}u-u^2\bigr)\mathcal{F}\tilde\mathcal{O}(u) \bigr),
\\
\tilde\varphi'_T(u-\mathrm{i})&:=&(\mathrm{i}-2u)\mathcal{F}\tilde\mathcal{O}(u)-\bigl(u+\mathrm{i}u^2\bigr)\mathcal F \bigl(x\tilde\mathcal{O}(x)
\bigr) (u),\qquad u\in\mathbb R.
\end{eqnarray*}
For all $u\in\mathbb R$ where $|\tilde\varphi_T(u-\mathrm{i})|>\kappa(u)$ we
obtain $\tilde\varphi_T(u-\mathrm{i})=1+(\mathrm{i}u-u^2)\mathcal{F}\tilde\mathcal{O}(u)$. For
$|\tilde\varphi_T(u-\mathrm{i})|=\kappa(u)$ the estimate $|\tilde\varphi
_T(u-\mathrm{i})-\varphi_T(u-\mathrm{i})|\geq2\kappa(u)$ follows from \eqref
{eqconditionKappa}. This yields
\begin{eqnarray*}
\bigl|\tilde\varphi_T(u-\mathrm{i})-\varphi_T(u-\mathrm{i})\bigr| &\leq&\bigl|1+
\bigl(\mathrm{i}u-u^2\bigr)\mathcal{F}\tilde\mathcal{O}(u)-\varphi_T(u-\mathrm{i})\bigr|+
\kappa(u)
\\
& \leq&\bigl|1+\bigl(\mathrm{i}u-u^2\bigr)\mathcal{F}\tilde\mathcal{O}(u)-
\varphi_T(u-\mathrm{i})\bigr|+\tfrac
{1}{2}\bigl|\tilde\varphi_T(u-\mathrm{i})-
\varphi_T(u-\mathrm{i})\bigr|.
\end{eqnarray*}
Therefore, $|\tilde\varphi_T(u-\mathrm{i})-\varphi_T(u-\mathrm{i})|\leq
2|1+(\mathrm{i}u-u^2)\mathcal{F}\tilde\mathcal{O}(u)-\varphi_T(u-\mathrm{i})|$ holds for all
$u\in
\mathbb R$. We obtain a similar decomposition as \cite{KappusReiss2010},
\begin{eqnarray*}
\bigl|\tilde\psi'(u)-\psi'(u)\bigr|& =&\frac{1}{T} \biggl|
\frac{\tilde\varphi'_T(u-\mathrm{i})}{\tilde\varphi_T(u-\mathrm{i})}-\frac{\varphi
'_T(u-\mathrm{i})}{\varphi_T(u-\mathrm{i})} \biggr|
\\
& \leq&\frac{1}{T|\tilde\varphi_T(u-\mathrm{i})|} \bigl(\bigl|\tilde\varphi'_T(u-\mathrm{i})-
\varphi'_T(u-\mathrm{i})\bigr|+T\bigl|\psi'(u)\bigr|\bigl|
\varphi_T(u-\mathrm{i})-\tilde\varphi_T(u-\mathrm{i})\bigr| \bigr)
\\
& \leq&\frac{1}{2T\kappa(u)} \bigl( \bigl(\bigl(1+4u^2\bigr)^{1/2}+2T\bigl|
\psi'(u)\bigr|\bigl(u^2+u^4\bigr)^{1/2}
\bigr)\bigl|\mathcal{F}(\tilde\mathcal{O}-\mathcal{O})(u)\bigr|
\\
& &\hspace*{38pt}{}+\bigl(u^2+u^4\bigr)^{1/2}\bigl|\mathcal F
\bigl(x(\tilde\mathcal{O}-\mathcal{O}) (x) \bigr) (u)\bigr| \bigr).
\end{eqnarray*}
Since $|\psi'(u)|\leq|\gamma|+\|k_e\|_{L^1}\leq2R$, we have
\[
\bigl|\tilde\psi'(u)-\psi'(u)\bigr|\lesssim\frac{1}{\kappa(u)}
\bigl(\bigl(1+u^2\bigr)\bigl|\mathcal{F}(\tilde\mathcal{O}-\mathcal{O})(u)\bigr| +
\bigl(u^2+u^4\bigr)^{1/2}\bigl|\mathcal F \bigl(x(
\tilde\mathcal{O}-\mathcal{O}) (x) \bigr) (u)\bigr| \bigr).
\]
It follows with Plancherel's equality
\begin{eqnarray*}
S&\leq&3\mathbb E \bigl[\bigl \|\mathcal F^{-1} \bigl(\bigl(\tilde
\psi'(u)-\psi'(u)\bigr)w_k(u/U) \bigr)
\bigr\|_{L^2}^2 \bigr] \\
&=&\frac{3}{2\uppi}\int
_{\mathbb R}\mathbb E \bigl[\bigl|\tilde\psi'(u)-
\psi'(u)\bigr|^2 \bigr]\bigl|w_k(u/U)\bigr|^2
\,\mathrm{d}u
\\
&\lesssim&\int_\R\frac{u^4}{|\kappa(u)|^2} \bigl(\mathbb E
\bigl[\bigl|\mathcal{F}(\tilde\mathcal{O}-\mathcal{O})(u)\bigr|^2 \bigr]\\
&&\hspace*{37pt}\quad{} +\mathbb E \bigl[ \bigl|
\mathcal F \bigl(x(\tilde\mathcal{O}-\mathcal{O}) (x) \bigr) (u) \bigr|^2 \bigr]
\bigr)\bigl|w_k(u/U)\bigr|^2\,\mathrm{d}u
\\
&=:&S_1+S_2.
\end{eqnarray*}
Both terms can be estimated similarly. Thus, we only write it down for
$S_2$, where stronger conditions are needed. Lemma \ref{lemxbBound}
and $\|\mathcal F(xb_j(x))\|_\infty\leq2\Delta(x_j+\Delta)$,
$j=1,\ldots
,N$, yield
\begin{eqnarray*}
S_2&\leq&\int_\R\frac{u^4}{|\kappa(u)|^2} \bigl(
\bigl\|x(O_l-O) (x)\bigr\|_\infty^2 +\var\bigl(\mathcal F
\bigl(x\tilde\mathcal{O}(x) \bigr) (u) \bigr) \bigr)\bigl|w_k(u/U)\bigr|^2
\,\mathrm{d}u
\\
&\lesssim&\int_\R|u|^{2T\bar\alpha+4} \Biggl(
\Delta^4+\sum_{j=1}^N
\delta_j^2\bigl|\mathcal F \bigl(xb_j(x) \bigr)
(u)\bigr|^2 \Biggr)\bigl|w_k(u/U)\bigr|^2\,\mathrm{d}u
\\
&\lesssim&\bigl(\Delta^4+\Delta^2\bigl\|(x_j
\delta_j)\bigr\|_{l^2}^2+\Delta^4\|
\delta_j\|_{l^2}^2\bigr)U^{2T\bar\alpha+5}
\lesssim\varepsilon^2U^{2T\bar\alpha+5}.
\end{eqnarray*}

Therefore, we have shown $\mathbb E [\|\hat k_e-k_e\|_{L^2,\tau
}^2 ]\lesssim U^{-2s}+\varepsilon^2U^{2T\bar\alpha+5}+U\mathbb
E[|\hat\gamma-\gamma|^2]$. The assertion follows from the asymptotic
optimal choice $U=U_{\bar\alpha}=\varepsilon^{-2/(2s+2T\bar\alpha+5)}$
and the assumption on the risk of $\hat\gamma$.

%s6.2 #&#
\subsection{\texorpdfstring{Proof of Proposition \protect\ref{propadaptiveEstimator}}
{Proof of Proposition 4.1}}

\textit{Step} 1: Let $(a_\varepsilon)_{\varepsilon>0}$ be a
deterministic sequence such that there is a constant $C>0$ with
$|a_\varepsilon-\alpha|\leq C|\log\varepsilon|^{-1}$. Let the
estimator $\hat\alpha_0$ use the cut-off value $ U_\varepsilon
:=\tilde
U_{a_\varepsilon}$ and the trimming parameter $\kappa_\varepsilon
:=\tilde\kappa_{\bar a_\varepsilon}$, with $\bar a_\varepsilon
:=a_\varepsilon+C|\log\varepsilon|^{-1}$, as defined in \eqref
{equEpsilon} and \eqref{eqkappaEpsilon}. Then we can show the
asymptotic risk bound $\sup_{\mathcal P\in\mathcal G_s(R,\alpha
)}\mathbb E_{\mathcal P}[|\hat\alpha_0-\alpha|^2]^{1/2} \lesssim
\varepsilon^{2(s-1)/(2s+2T\alpha+1)}$
as follows:
By construction holds $\alpha\leq\bar a_\varepsilon$. Hence, $\kappa
_\varepsilon$ fulfills condition \eqref{eqconditionKappa} for each
pair $\mathcal P\in\mathcal G_s(R,\alpha)$ and thus we deduce from
Theorem \ref{thmupperBoundsParameters}
%
%
%e6.1 #&#
\begin{eqnarray}\label{eqinDetermRes}
&&\mathbb E_{\mathcal P} \bigl[\bigl|\hat\alpha_0-
\alpha\bigr|^2 \bigr]\nonumber\\
 && \quad\lesssim U_\varepsilon^{-2(s-1)}+U_\varepsilon
^{2T\alpha
+3}
\varepsilon^2 +U_\varepsilon^{4T\bar a_\varepsilon+8}\varepsilon^4
\\
&&\quad=\varepsilon^{4(s-1)/(2s+2Ta_\varepsilon+1)} \bigl(1+ \varepsilon
^{4T(a_\varepsilon-\alpha)/(2s+2Ta_\varepsilon+1)} +
\varepsilon^{(4s-8+8T(a_\varepsilon-\bar a_\varepsilon
))/(2s+2Ta_\varepsilon+1)} \bigr).\nonumber
\end{eqnarray}
The first factor has the claimed order, since $\varepsilon
^{4(s-1)/(2s+2Ta_\varepsilon+1)}\lesssim\varepsilon^{4(s-1)/(2s+2T\alpha
+1)}$ follows with easy calculations from $(\alpha-a_\varepsilon)\log
\varepsilon\leq C$.
% \begin{eqnarray*}\begin{array}{rrl}
% &(\alpha-a_\varepsilon)\log\varepsilon&\leq C\\
% \Rightarrow& (2s+2T\alpha+1)\log\varepsilon&\leq(2s+2Ta_
% \Rightarrow& \frac{4(s-1)}{2s+2Ta_\varepsilon+1}\log\varepsilon&\leq
% \Rightarrow& \varepsilon^{4(s-1)/(2s+2Ta_\varepsilon+1)} &\lesssim
% \end{array}
% \end{eqnarray*}
Hence, the claim follows once we have bound the sum in the bracket of
equation \eqref{eqinDetermRes}. For the second term, this is implied by
\[
\biggl|\frac{4T(a_\varepsilon-\alpha)}{2s+2Ta_\varepsilon+1}\log\varepsilon
\biggr| \leq\frac{4T|(a_\varepsilon-\alpha)\log\varepsilon|}{2s+1} \leq
\frac{4TC}{2s+1}.
\]
To estimate the third term, we obtain from $s\geq2$ and $\varepsilon<1$
\[
\frac{4s-8+8T(a_\varepsilon-\bar a_\varepsilon)}{2s+2Ta_\varepsilon
+1}\log\varepsilon\leq\frac{-8TC|\log\varepsilon|^{-1}}{2s+1}\log
\varepsilon\leq
\frac
{8TC}{2s+1}.
\]
\textit{Step} 2: Let $\mathcal P\in\mathcal G_s(R,\alpha)$. Note that
$\kappa_\varepsilon$ satisfies the condition \eqref{eqconditionKappa}
on the set $\{|\hat\alpha_{\mathrm{pre}}-\alpha|<|\log\varepsilon|^{-1}\}$.
Using the independence of $\hat\alpha_{\mathrm{pre}}$ and $O_j$, the almost sure
bound $\tilde\alpha_0\leq\bar\alpha$ and the concentration of
$\hat
\alpha_{\mathrm{pre}}$, we deduce from step 1:
\begin{eqnarray*}
\mathbb E_{\mathcal P,\hat\alpha_{\mathrm{pre}}} \bigl[\bigl|\tilde\alpha_0-\alpha
\bigr|^2 \bigr] &\leq&\mathbb E_{\mathcal P,\hat\alpha_{\mathrm{pre}}} \bigl[\mathbb
E_{\mathcal
P,\hat\alpha_{\mathrm{pre}}} \bigl[|\tilde\alpha_0-\alpha|^2 |\hat
\alpha_{\mathrm{pre}} \bigr]\mathbf1_{\{|\hat\alpha_{\mathrm{pre}}-\alpha|<|\log
\varepsilon
|^{-1}\}} \bigr]
\\
&&{} +4\bar\alpha^2\mathbb P_{\hat\alpha_{\mathrm{pre}}} \bigl(|\hat
\alpha_{\mathrm{pre}}-\alpha|\geq|\log\varepsilon|^{-1} \bigr)
\\
&\lesssim& \varepsilon^{4(s-1)/(2s+2T\alpha+1)}+4\bar\alpha^2{d}
\varepsilon^2.
\end{eqnarray*}
Since the second term decreases faster then the first one for
$\varepsilon\to0$, we obtain the claimed rate.

%s6.3 #&#
\subsection{\texorpdfstring{Proof of Proposition \protect\ref{propconcentrationAlpha}}{Proof of Proposition 4.2}}

Recall that the cut-off value of $\hat\alpha_0$ is given by
$U=\varepsilon^{-2/(2s+2T\bar\alpha+1)}$. For $\kappa>0,$ we obtain
from the definition of the estimator and the decomposition of the
stochastic error into linear part and remainder:
\begin{eqnarray*}
\mathbb P\bigl(|\hat\alpha_0-\alpha|\geq\kappa\bigr) &=&\mathbb P \biggl( \biggl|
\int_{-U}^U\re(\rho+\tilde\psi-\psi)
(u)w_{\alpha
_0}^U(u)\,\mathrm{d}u \biggr|\geq\kappa\biggr)
\\
&\leq&\mathbb P \biggl( \biggl|\int_{-U}^U
\rho(u)w_{\alpha_0}^U(u)\, \mathrm{d}u \biggr|\geq\frac{\kappa}{3}
\biggr) +\mathbb P \biggl( \biggl|\int_{-U}^U\re\bigl(
\mathcal L(u)\bigr)w_{\alpha
_0}^U(u)\,\mathrm{d} u \biggr|\geq
\frac{\kappa}{3} \biggr)
\\
&&{} +\mathbb P \biggl( \biggl|\int_{-U}^U\mathcal
R(u)w_{\alpha_0}^U(u)\, \mathrm{d} u \biggr|\geq\frac{\kappa}{3}
\biggr)=:P_1+P_2+P_3.
\end{eqnarray*}
We will bound all three probabilities separately. To that end, let
$c_j,j\in\mathbb N,$ be suitable non-negative constants not depending
on $\kappa, \varepsilon$ and $N$.

The event in $P_1$ is deterministic. Hence, the same estimate on the
deterministic error as in Theorem \ref{thmupperBoundsParameters}
\[
\biggl|\int_{-U}^U\rho(u)w_{\alpha_0}^U(u)
\,\mathrm{d}u \biggr|\leq c_1U^{-(s-1)}=c_1
\varepsilon^{2(s-1)/(2s+2T\bar\alpha+1)}
\]
yields $P_1=0$ for all $\varepsilon<\varepsilon^{(1)}:= (\kappa
/(3c_1) )^{(2s+2T\bar\alpha+1)/(2s-2)}$.

To bound $P_2$ we infer from the definition of $\mathcal L$, the
linearity of the errors in $\tilde\mathcal{O}=\mathcal O_l+\sum
_{j=1}^N\delta_j\varepsilon_jbj$ and from the estimate of the term
$|\mathcal L_b|$ in Theorem \ref{thmupperBoundsParameters}
\begin{eqnarray*}
&& \biggl|\int_{-U}^U\re\bigl(\mathcal L(u)
\bigr)w_{\alpha_0}^U(u)\,\mathrm{d}u \biggr|\\
&&\quad = \biggl|\int
_{-U}^U\re\biggl(\frac{(\mathrm{i}-u)u}{T\varphi_T(u-\mathrm{i})}\mathcal{F}(
\tilde\mathcal{O}-\mathcal O) (u) \biggr)w_{\alpha_0}^U(u)\,\mathrm{d}u
\biggr|
\\
&&\quad \leq\int_{-U}^U\frac{(u^4+u^2)^{1/2}}{T|\varphi_T(u-\mathrm{i})|}\bigl|\mathcal
{F(O}_l-\mathcal O) (u)w_{\alpha_0}^U(u)\bigr|\,
\mathrm{d}u
\\
&&\qquad{} + \biggl|\int_{-U}^U\re\Biggl(\frac{(\mathrm{i}-u)u}{T\varphi_T(u-\mathrm{i})}
\sum_{j=1}^N\delta_j
\varepsilon_j\mathcal Fb_j(u) \Biggr)w_{\alpha
_0}^U(u)
\,\mathrm{d} u \biggr|
\\
&&\quad \leq c_2\Delta^2U^{T\bar\alpha+2} + \Biggl|\sum
_{j=1}^N\delta_j\varepsilon_j
\int_{-U}^U\re\biggl(\frac{(\mathrm{i}-u)u}{T\varphi_T(u-\mathrm{i})}\mathcal
Fb_j(u) \biggr)w_{\alpha_0}^U(u)\,\mathrm{d}u \Biggr|
\\
&&\quad \leq c_2\varepsilon^{2(s-1)/(2s+2T\bar\alpha+1)} + \Biggl|\sum
_{j=1}^Na_j\varepsilon_j \Biggr|,
\end{eqnarray*}
where the coefficients are given by $a_j:=\delta_j\int_{-U}^U\re
(\frac{(\mathrm{i}-u)u}{T\varphi_T(u-\mathrm{i})}\mathcal Fb_j(u) )w_{\alpha
_0}^U(u)\,\mathrm{d}u$ for $j=1,\ldots,N$. To apply \eqref
{eqconcentrationEpsilon}, we deduce from $\|\mathcal Fb_j\|_\infty
\leq
2\Delta$, the weight function property \eqref
{eqweightfunctionProperty} and the assumption $\Delta\|\delta\|
_{l^2}^2\lesssim\|\delta\|_{l^\infty}^2$
\begin{eqnarray*}
\sum_{j=1}^Na_j^2
&\leq&\sum_{j=1}^N\delta_j^2
\biggl(\int_{-U}^U\frac
{(u^4+u^2)^{1/2}}{T|\varphi_T(u-\mathrm{i})|}\bigl|\mathcal
Fb_j(u)\bigr|\bigl|w_{\alpha
_0}^U(u)\bigr|\,\mathrm{d}u
\biggr)^2 \leq c_3\Delta^2U^{2T\bar\alpha+4}\|
\delta\|_{l^2}^2
\\
&\leq& c_4\varepsilon^2U^{2T\bar\alpha+4} =c_4
\varepsilon^{2(s-1)/(2s+2T\bar\alpha+1)}.
\end{eqnarray*}
This implies through the concentration inequality of $(\varepsilon_j)$
\begin{eqnarray*}
P_2&\leq&\mathbb P \Biggl(\Biggl |\sum_{j=1}^Na_j
\varepsilon_j \Biggr|\geq\frac
{\kappa}{6} \Biggr) +\mathbb P
\biggl(c_2\varepsilon^{2(s-1)/(2s+2T\bar
\alpha
+1)}\geq\frac{\kappa}{6} \biggr)
\\
&\leq&C_1\exp\biggl(-\frac{C_2}{36c_4}\kappa^2
\varepsilon^{-2(s-1)/(2s+2T\bar\alpha+1)} \biggr)
\end{eqnarray*}
for all $\varepsilon<\varepsilon^{(2)}:= (\kappa/(6c_2)
)^{(2s+2T\bar\alpha+1)/(2s-2)}$.

It remains to estimate probability $P_3$. The bound of $\mathcal R$ in
Proposition \ref{propresultsReiss} ii) yields
\begin{eqnarray*}
 &&\biggl|\int_{-U}^U\mathcal R(u)w_{\alpha_0}^U(u)
\,\mathrm{d}u \biggr| \\
&&\quad\leq\int_{-U}^U
\frac{u^4+u^2}{T\kappa(u)^2}\bigl|\mathcal{F}(\tilde\mathcal{O}-\mathcal{O})(u)\bigr|^2\bigl|w_{\alpha_0}^U(u)\bigr|
\,\mathrm{d}u
\\
&&\quad\leq2\int_{-U}^U\frac{u^4+u^2}{T\kappa(u)^2}\bigl|\mathcal
{F(O}_l-\mathcal O) (u)\bigr|^2\bigl|w_{\alpha_0}^U(u)\bigr|
\,\mathrm{d}u \\
&&\qquad{}+2\int_{-U}^U\frac{u^4+u^2}{T\kappa(u)^2} \Biggl|
\sum_{j=1}^N\delta_j
\varepsilon_j\mathcal Fb_j(u)\Biggr |^2\bigl|w_{\alpha_0}^U(u)\bigr|
\,\mathrm{d}u.
\end{eqnarray*}
The first addend gets small owing to Proposition \ref
{propresultsReiss}(i):
\begin{eqnarray*}
&&\int_{-U}^U\frac{u^4+u^2}{T\kappa(u)^2}\bigl|
\mathcal{F(O}_l-\mathcal O) (u)\bigr|^2\bigl|w_{\alpha_0}^U(u)\bigr|
\,\mathrm{d}u \\
&&\quad\leq\bigl\|\mathcal{F(O}_l-\mathcal O)\bigr\|_\infty^2
\int_{-U}^U\frac
{u^4+u^2}{T\kappa(u)^2}\bigl|w_{\alpha_0}^U(u)\bigr|
\,\mathrm{d}u
\\
&&\quad\leq c_5\Delta^4U^{2T\bar\alpha+4} \leq
c_5\varepsilon^{2(s-1)/(2s+2T\bar\alpha+1)}.
\end{eqnarray*}
For the second one, we obtain
\begin{eqnarray*}
\Biggl|\sum_{j=1}^N\delta_j
\varepsilon_j\mathcal Fb_j(u) \Biggr|^2 =\sum
_{j=1}^N\delta_j^2
\varepsilon_j^2\bigl|\mathcal Fb_j(u)\bigr|^2
+2\sum_{j=2}^N\sum
_{k=1}^{j-1}\delta_j\delta_k
\varepsilon_j\varepsilon_k\re\bigl(\mathcal
Fb_j(u)\mathcal Fb_k(-u) \bigr).
\end{eqnarray*}
Thus,
\begin{eqnarray*}
\biggl|\int_{-U}^U\mathcal R(u)w_{\alpha_0}^U(u)
\,\mathrm{d}u \biggr| \leq2c_5\varepsilon^{{(4s-6)}/{(2s+2T\bar\alpha+1)}}
+2\sum
_{j=1}^N\delta_j^2
\varepsilon_j^2\xi_{j,j}(U)+4\sum
_{j=2}^N\sum_{k=1}^{j-1}
\delta_j\delta_k\varepsilon_j
\varepsilon_k\xi_{j,k}(U)
\end{eqnarray*}
with $\xi_{j,k}(U):=\int_{-U}^U\frac{u^4+u^2}{T\kappa(u)^2}\re
(\mathcal Fb_j(u)\mathcal Fb_k(-u) )|w_{\alpha_0}^U(u)|\,\mathrm{d}u$.
Denoting the diagonal term and the cross term as
\[
D_N:=\sum_{j=1}^N
\delta_j^2\varepsilon_j^2
\xi_{j,j}(U)\quad \mbox{and}\quad U_N:=\sum
_{j=2}^N\sum_{k=1}^{j-1}
\delta_j\delta_k\varepsilon_j
\varepsilon_k\xi_{j,k}(U),
\]
respectively, we obtain
\[
P_3\leq\mathbb P \biggl(2c_5\varepsilon^{2(s-1)/(2s+2T\bar\alpha
+1)}
\geq\frac{\kappa}{9} \biggr) +\mathbb P \biggl(2D_N\geq
\frac{\kappa}{9} \biggr)+\mathbb P \biggl(4U_N\geq
\frac{\kappa}{9} \biggr).
\]
The first summand vanishes for $\varepsilon<\varepsilon^{(3)}:=
(\kappa/(18c_5) )^{(2s+2T\bar\alpha+1)/(2s-2)}$. To estimate the
probabilities on $D_N$ and $U_N$, we establish the bound
%
%
%e6.2 #&#
\begin{equation}
\label{eqinConcentrationXi}\bigl |\xi_{j,k}(U)\bigr| \leq\|\mathcal
Fb_j\|_\infty\|\mathcal Fb_k\|_\infty
\int_{-U}^U\frac
{u^4+u^2}{T\kappa(u)^2}\bigl|w_{\alpha_0}^U(u)\bigr|
\,\mathrm{d}u \leq c_6\Delta^2U^{2T\bar\alpha+4}
\end{equation}
for $j,k=1,\ldots,N$. Hence,
\[
\Biggl|\sum_{j=1}^N\delta_j^2
\xi_{j,j}(U) \Biggr| \leq c_6\Delta^2\|\delta
\|_{l^2}^2U^{2T\bar\alpha+4} \leq c_7
\varepsilon^2U^{2T\bar\alpha+4} \leq c_7
\varepsilon^{2(s-1)/(2s+2T\bar\alpha+1)},
\]
which yields together with \eqref{eqconcentrationEpsilon}
\begin{eqnarray*}
\mathbb P \biggl(D_N\geq\frac{\kappa}{18} \biggr) &\leq&\mathbb P
\Biggl(\sup_{k=1,\ldots,N}|\varepsilon_k|^2 \Biggl|\sum
_{j=1}^N\delta_j^2
\xi_{j,j}(U) \Biggr|\geq\frac{\kappa}{18} \Biggr)
\\
&\leq&\mathbb P \biggl(\sup_{k=1,\ldots,N}|\varepsilon_k|^2
\geq\frac
{\kappa
}{18c_7}\varepsilon^{-2(s-1)/(2s+2T\bar\alpha+1)} \biggr)
\\
&\leq& C_1N\exp\biggl(-\frac{C_2}{18c_7}\kappa\varepsilon
^{-2(s-1)/(2s+2T\bar\alpha+1)}
\biggr).
\end{eqnarray*}
To derive an exponential inequality for the U-statistic $U_N$, we apply
the martingale idea in \cite{houdreReynaud2003}. Because of the
independence and the centering of the $(\varepsilon_j)$, the process
$(U_N)_{N\geq1}$ is a martingale with respect to its natural filtration
$(\mathcal F_N^U)$ (setting $U_1=0$):
\[
\mathbb E\bigl[U_N-U_{N-1}|\mathcal F_{N-1}^U
\bigr] =\mathbb E \Biggl[\sum_{k=1}^{N-1}
\delta_N\delta_k\varepsilon_N
\varepsilon_k\xi_{N,k}(U)\Big|\mathcal F_{N-1}^U
\Biggr]=0.
\]
We apply the martingale version of the Bernstein
inequality, see Theorem~VII.3.6 in \cite{shiryaev1995}, which yields for arbitrary $t,Q,S>0$
%
%
%e6.3 #&#
\begin{eqnarray}
\label{eqBernstein} \mathbb P\bigl(|U_N|\geq t\bigr) &\leq&2\mathbb P\bigl(
\anglel U\angler_N>Q\bigr)+2\mathbb P \Bigl(\max_{k=1,\ldots
,N}|U_k-U_{k-1}|>S
\Bigr)
\nonumber
\\[-8pt]
\\[-8pt]
\nonumber
 &&{}+2\exp\biggl(-\frac{t^2}{4(Q+t S)} \biggr).
\end{eqnarray}
Hence, we consider the increment $|U_N-U_{N-1}|=|\varepsilon_N|
|\sum_{k=1}^{N-1}\delta_N\delta_k\xi_{N,k}(U)\varepsilon_k |$, for
$N\geq
2$. Denoting $a_{N,k}:=\delta_N\delta_k\xi_{N,k}(U)$, we estimate using
\eqref{eqinConcentrationXi}
%
%
%e6.4 #&#
\begin{eqnarray}\label{eqinConcentrationANk}
\sum_{k=1}^{N-1}a_{N,k}^2
&=&\delta_N^2\sum_{k=1}^{N-1}
\delta_k^2\xi_{N,k}(U)^2 \leq
c_6^2\Delta^4U^{4T\bar\alpha+8}
\delta_N^2\|\delta\|_{l^2}^2
\nonumber
\\[-8pt]
\\[-8pt]
\nonumber
&\leq& c_6^2\Delta^4\|\delta
\|_{l^2}^4U^{4T\bar\alpha+8} \leq c_7^2
\varepsilon^4U^{4T\bar\alpha+8} \leq c_7^2
\varepsilon^{4(s-1)/(2s+2T\bar\alpha+1)}
.
\end{eqnarray}
Thus, by Assumption \eqref{eqconcentrationEpsilon} we obtain for all $S>0$
\begin{eqnarray*}
&&\mathbb P\bigl(|U_N-U_{N-1}|> S\bigr)\\
&&\quad=\mathbb P \Biggl(|
\varepsilon_N| \Biggl|\sum_{k=1}^{N-1}a_{N,k}
\varepsilon_k\Biggr |>S \Biggr)
\\
&&\quad\leq\mathbb P \bigl(|\varepsilon_N|> \sqrt{S}\varepsilon
^{-(s-1)/(2s+2T\bar\alpha+1)}
\bigr) +\mathbb P \Biggl( \Biggl|\sum_{k=1}^{N-1}a_{N,k}
\varepsilon_k\Biggr |> \sqrt{S}\varepsilon^{(s-1)/(2s+2T\bar\alpha+1)} \Biggr)
\\
&&\quad\leq C_1\exp\bigl(-C_2S\varepsilon^{-2(s-1)/(2s+2T\bar\alpha
+1)}
\bigr)\\
&&\qquad{} +C_1\exp\biggl(-\frac{C_2}{c_7^2}S\varepsilon
^{-2(s-1)/(2s+2T\bar
\alpha
+1)}
\biggr).
\end{eqnarray*}
The quadratic variation of $U_N$ is given by
\begin{eqnarray*}
\anglel U\angler_N-\anglel U\angler_{N-1} =\mathbb E
\bigl[(U_N-U_{N-1})^2 |\mathcal
F_{N-1}^U \bigr] =\delta_N^2
\Biggl(\sum_{k=1}^{N-1}\delta_k
\varepsilon_k\xi_{N,k}(U) \Biggr)^2.
\end{eqnarray*}
W.l.o.g. we can assume $\sum_{j=2}^N\delta_j^2>0$. Otherwise follows
$\sum_{j=2}^N\delta_j^2=0$ which implies $\delta_j=0$ for all
$j=2,\ldots
,N$ and thus $\anglel U\angler_N=\sum_{j=2}^N (\anglel U\angler
_j-\anglel U\angler_{j-1} )=0$. Then $\mathbb P(\anglel U\angler
_N>Q)=0$ would
hold for $Q>0$. Hence, we obtain:
\begin{eqnarray*}
\mathbb P\bigl(\anglel U\angler_N>Q\bigr) &=&\mathbb P \Biggl(\sum
_{j=2}^N \bigl(\anglel U
\angler_j-\anglel U\angler_{j-1} \bigr)>Q \Biggr) \leq\sum
_{j=2}^N\mathbb P \biggl(\anglel U
\angler_j-\anglel U\angler_{j-1}>\frac{\delta_j^2}{\sum_{k=2}^N\delta
_k^2}Q \biggr)
\\
% &\leq\sum_{j=2}^N\mathbb P\Big(\delta_j^2\Big(\sum_{k=1}^{j-1}\delta_k
&\leq&\sum_{j=2}^N
\mathbb P \Biggl(\|\delta\|_{l^2}\sum_{k=1}^{j-1}
\delta_k\varepsilon_k\xi_{j,k}(U)>\sqrt{Q}
\Biggr).
\end{eqnarray*}
To apply inequality \eqref{eqconcentrationEpsilon}, we estimate $\|
\delta\|_{l^2}^2\sum_{k=1}^{j-1}\delta_k^2\xi_{j,k}(U)^2\leq
c_6^2\Delta^4\|\delta\|_{l^2}^4U^{4T\bar\alpha+8}\leq  c_7^2\varepsilon
^{4(s-1)/(2s+2T\bar\alpha+1)}$ analogous to \eqref
{eqinConcentrationANk} and obtain
\[
\mathbb P\bigl(\anglel U\angler_N>Q\bigr)\leq C_1N\exp
\biggl(-\frac
{C_2}{c_7^2}Q\varepsilon^{-4(s-1)/(2s+2T\bar\alpha+1)} \biggr).
\]
We deduce from Bernstein's inequality \eqref{eqBernstein}
\begin{eqnarray*}
&&\mathbb P \biggl(U_N\geq\frac{\kappa}{36} \biggr)
\\
&&\quad\leq 2\mathbb P \bigl(\anglel U\angler_N>Q \bigr)+2\mathbb P
\Bigl(\max_{k=2,\ldots,N}|U_k-U_{k-1}|>S \Bigr)+2\exp
\biggl(-\frac{\kappa^2}{144(36Q+\kappa S)} \biggr)
\\
&&\quad\leq 2C_1N\exp\biggl(-\frac{C_2}{c_7^2}Q\varepsilon
^{-4(s-1)/(2s+2T\bar
\alpha+1)}
\biggr)
\\
&&\qquad{}+ 4C_1N\exp\biggl(-\frac{C_2}{c_7^2\vee1}S\varepsilon
^{-2(s-1)/(2s+2T\bar\alpha+1)}
\biggr) + 2\exp\biggl(-\frac{\kappa^2}{144(36Q+\kappa S)} \biggr).
\end{eqnarray*}
By choosing $Q=\kappa S$ and $S=\sqrt{\kappa}\varepsilon
^{(s-1)/(2s+2T\bar\alpha+1)}$, we get
\[
\mathbb P \biggl(U_N\geq\frac{\kappa}{36} \biggr)
\leq(6C_1N+2)\exp\Bigl(-c_8\min_{q=1,3} \bigl(
\kappa^{1/2}\varepsilon^{-(s-1)/(2s+2T\bar\alpha+1)} \bigr)^q \Bigr).
\]
For all $\varepsilon<\varepsilon^{(3)}$, we have $\kappa\varepsilon
^{-2(s-1)/(2s+2T\bar\alpha+1)}> \kappa(\varepsilon^{(3)}
)^{-2(s-1)/(2s+2T\bar\alpha+1)}\sim1$ and hence,
\begin{eqnarray*}
P_3\leq\mathbb P\biggl(D_N\geq\frac{\kappa}{18}
\biggr)+\mathbb P\biggl(U_N\geq\frac
{\kappa}{36}\biggr)
\leq(7C_1N+2)\exp\bigl(-c_8\kappa^{1/2}
\varepsilon^{-(s-1)/(2s+2T\bar
\alpha+1)} \bigr).
\end{eqnarray*}
Putting the bounds of $P_1, P_2$ and $P_3$ together yields for a
constant $c\in(0,\infty)$ and all $\varepsilon<\varepsilon_0\wedge1$
with $\varepsilon_0:=\min\{\varepsilon^{(1)}, \varepsilon
^{(2)},\varepsilon^{(3)}\}$
\[
\mathbb P\bigl(|\hat\alpha_0-\alpha|\geq\kappa\bigr) \leq(7C_1N+C_1+2)
\exp\bigl(-c\bigl(\kappa^2\wedge\kappa^{1/2}\bigr)
\varepsilon^{-(s-1)/(2s+2T\bar\alpha+1)} \bigr).
\]

%sA #&#
\begin{appendix}
\section*{\texorpdfstring{Appendix: Proof of Lemma \protect\ref{lemboundPhi}}
{Appendix: Proof of Lemma 2.1}}\label{app}
\textit{Part} (i) The martingale condition yields
\begin{eqnarray*}
\bigl|\varphi_T(u-\mathrm{i})\bigr|=\exp\biggl(T\int_{-\infty}^\infty
\bigl(\cos(ux)-1 \bigr)\frac{\mathrm{e}^xk(x)}{|x|}\,\mathrm{d}x \biggr).
\end{eqnarray*}
W.l.o.g. we assume $T=1$, $\alpha>0$ and $u\geq1$ because of the
symmetry of the cosine.

We split the integral domain into three parts:
%
%
%eA.1 #&#
\setcounter{equation}{0}
\begin{equation}
\label{eqApp} \bigl|\varphi_1(u-\mathrm{i})\bigr|=\exp\biggl( \biggl(\int
_0^1+\int_1^{u}+
\int_{u}^\infty\biggr)\frac{\cos x-1}{x}
\biggl(\mathrm{e}^{x/u}k\biggl(\frac{x}{u}\biggr)+\mathrm{e}^{-x/u}k
\biggl(-\frac
{x}{u}\biggr) \biggr)\,\mathrm{d}x \biggr).
\end{equation}
Using $\|k\|_\infty\leq\|k\|_{\mathrm{TV}}<\infty$ by assumption and the
constant $C_1:=\int_0^1\frac{1-\cos x}{x}\,\mathrm{d}x\in(0,\infty
)$, we estimate
\[
\int_0^1\frac{\cos x-1}{x}
\biggl(\mathrm{e}^{x/u}k\biggl(\frac
{x}{u}\biggr)+\mathrm{e}^{-x/u}k
\biggl(-\frac
{x}{u}\biggr) \biggr)\,\mathrm{d}x\geq2\mathrm{e}^{1/u}\|k
\|_\infty\int_0^1\frac
{\cos
x-1}{x}\,
\mathrm{d}x \geq-2C_1e\|k\|_\infty.
\]
In the second part the dependence on $u$ comes into play. Writing
$\tilde k(x):=k(x)+k(-x)$, the Taylor series of the exponential
function together with dominated convergence yield
\begin{eqnarray*}
&&\int_1^{u}\frac{\cos x-1}{x}
\biggl(\mathrm{e}^{x/u}k\biggl(\frac
{x}{u}\biggr)+\mathrm{e}^{-x/u}k
\biggl(-\frac{x}{u}\biggr) \biggr)\,\mathrm{d}x
\\
&&\quad\geq\int_1^{u}\frac{\cos x-1}{x}\tilde k
\biggl(\frac{x}{u}\biggr)\,\mathrm{d}x +\|k\|_\infty\sum
_{k=1}^\infty\int_1^u(
\cos x-1)\frac
{x^{k-1}}{u^kk!}\,\mathrm{d}x
\\
&&\quad\geq-\alpha\log(u)+\int_{1/u}^{1} \bigl(\alpha-
\tilde k(x) \bigr)\frac
{\,\mathrm{d}x}{x}+\int_{1/u}^1
\frac{\cos(ux)}{x}\tilde k(x)\, \mathrm{d}x-2\|k\|_\infty\sum
_{k=1}^\infty\frac{1}{k!k}\bigl(1-u^{-k}
\bigr)
\\
&&\quad\geq-\alpha\log(u)-\sup_{x\in(0,1]}\frac{\tilde k(x)-\alpha
}{x}-2e\| k
\|_\infty+\int_{1/u}^1\frac{\cos(ux)}{x}
\tilde k(x)\,\mathrm{d}x.
\end{eqnarray*}
To bound the last term in the above display, we proceed as
Lemma 53.9 in \cite{sato1999}. By the bounded variation of $k$, we can define a
bounded signed measure $\rho$ via $\rho((a,b])=\tilde k(b+)-\tilde
k(a+), 0\le a<b$. Noting that $\int_y^\infty\frac{\cos x}{x}\,\mathrm
{d}x$ can
be bounded uniformly $y\in[1,\infty)$\vspace*{2pt} with a constant $C_2>0$, Fubini's
theorem yields
\begin{eqnarray*}
\int_{1/u}^1\frac{\cos(ux)}{x}\tilde k(x)\,
\mathrm{d}x& =&\tilde k\biggl(\frac{1}{u}+\biggr)\int_1^u
\frac{\cos x}{x}\,\mathrm{d}x+\int_{1}^u
\frac
{\cos x}{x}\int_{1/u}^{x/u}\rho(\mathrm{d}y)
\,\mathrm{d}x
\\
&=&\tilde k\biggl(\frac{1}{u}+\biggr)\int_1^u
\frac{\cos x}{x}\,\mathrm{d}x+\int_{1/u}^1\int
_{uy}^{u}\frac{\cos x}{x}\,\mathrm{d}x\,\rho(
\mathrm{d}y)
\\
&\geq&\tilde k\biggl(\frac{1}{u}+\biggr)\underbrace{\min_{v\ge1}
\int_1^v\frac
{\cos
x}{x}\,
\mathrm{d}x}_{=:-C_3\leq0} -2C_2\int_0^1|
\rho|(\mathrm{d}y) \\
&\geq&-2C_3\|k\|_\infty-4C_2
\|k\|_{\mathrm{TV}}.
\end{eqnarray*}
Obtaining for the third part in \eqref{eqApp} $\int_1^\infty\frac
{\cos
(ux)-1}{x} (\mathrm{e}^xk(x)+\mathrm{e}^{-x}k(x) )\,\mathrm{d}x\geq-2\|\mathrm{e}^xk(x)\|_{L^1}$, we
have with $q_k$ as defined in Lemma \ref{lemboundPhi}
\[
\bigl|\varphi_1(u-\mathrm{i})\bigr|\geq\exp\bigl(-q_k-(2e+4C_2+2C_3)
\|k\|_{\mathrm{TV}}-2\bigl\| \mathrm{e}^xk(x)\bigr\|_{L^1}
\bigr)u^{-\alpha}.
\]
We deduce the estimate $|\varphi_T(u-\mathrm{i})|\geq C_\varphi(T,q_k,\|
\mathrm{e}^xk(x)\|_{L^1},\|k\|_{\mathrm{TV}} \} )|u|^{-T\alpha}$ for $|u|\geq1$
with $C_\varphi(T,R):=\exp(-TR(3+2e+4C_2+2C_3))$.

\textit{Part} (ii) follows immediately from the explicit choice of
$C_\varphi$.
\end{appendix}

\section*{Acknowledgements}
The author thanks Markus Rei{\ss} and Jakob S\"ohl for providing many
helpful ideas and comments. The research was supported by the
Collaborative Research Center 649 ``Economic Risk'' of the German
Research Foundation (Deutsche Forschungsgemeinschaft).

\begin{supplement}%[id=suppA]
\stitle{Characteristic exponent and lower risk bounds}
\slink[doi]{10.3150/12-BEJ478SUPP} %[doi,text={...}] - jei reikia
%suskaldyti doi
\sdatatype{.pdf}
\sfilename{BEJ478\_supp.pdf}
\sdescription{First, we derive the representation of the characteristic
exponent given in Proposition~\ref{proprepPsi}. Furthermore, we
discuss Le Cam's asymptotic equivalence of our nonparametric regression
model to the continuous-time white noise model and show lower bounds in
the latter one.}
\end{supplement}

% imsref loaded by akundreckaite, 2012-12-13 15:59:23
%
% imsref loaded by akundreckaite, 2012-12-14 08:55:13
%
% imsref loaded by akundreckaite, 2012-12-14 10:31:22
% imsref loaded by akundreckaite, 2012-12-14 10:32:06
% imsref loaded by akundreckaite, 2012-12-17 11:09:39

%

\printhistory

\end{document}